\pdfoutput=1
\documentclass[tightenlines, twocolumn, notitlepage, nofootinbib, superscriptaddress]{revtex4}

\usepackage[dvipsnames]{xcolor}
\usepackage[intlimits]{amsmath}
\usepackage{amssymb}
\usepackage{bm}
\usepackage{bibunits}
\usepackage{enumerate}
\usepackage{floatrow}
\usepackage{hyperref}
\usepackage[symbol]{footmisc}
\usepackage{graphicx}
\usepackage{units}
\usepackage{xspace}
\usepackage{listings}
\usepackage{makecell}
\usepackage[utf8x]{inputenc}

\usepackage{floatrow}
\newfloatcommand{capbtabbox}{table}[][\FBwidth]
\usepackage{blindtext}

\DeclareMathOperator*{\argmin}{argmin}

\graphicspath{{./}, {figures/}}
\setcitestyle{super}

\newcommand{\golem}{\textsc{Golem}\xspace}

\makeatletter
\renewcommand*{\p@subsection}{\thesection.}
\makeatother

\makeatletter
\def\maketitle{
\@author@finish
\title@column\titleblock@produce
\suppressfloats[t]}
\makeatother

\begin{document}


	\title{\large{Golem: An algorithm for robust experiment and process optimization}}


    	\author{Matteo Aldeghi}
    	\email{matteo.aldeghi@vectorinstitute.ai}
	\affiliation{Vector Institute for Artificial Intelligence, Toronto, ON, Canada}
	\affiliation{Chemical Physics Theory Group, Department of Chemistry, University of Toronto, Toronto, ON, Canada}
	\affiliation{Department of Computer Science, University of Toronto, Toronto, ON, Canada}
	\author{Florian H\"ase}
	\affiliation{Department of Chemistry and Chemical Biology, Harvard University, Cambridge, MA, USA}
	\affiliation{Vector Institute for Artificial Intelligence, Toronto, ON, Canada}
	\affiliation{Chemical Physics Theory Group, Department of Chemistry, University of Toronto, Toronto, ON, Canada}
	\affiliation{Department of Computer Science, University of Toronto, Toronto, ON, Canada}
	\author{Riley J. Hickman}
	\affiliation{Chemical Physics Theory Group, Department of Chemistry, University of Toronto, Toronto, ON, Canada}
	\affiliation{Department of Computer Science, University of Toronto, Toronto, ON, Canada}
	\author{Isaac Tamblyn}
	\affiliation{National Research Council of Canada, Ottawa, ON, Canada}
	\affiliation{Vector Institute for Artificial Intelligence, Toronto, ON, Canada}
	\author{Al\'an Aspuru-Guzik}
	\email{alan@aspuru.com}
	\affiliation{Vector Institute for Artificial Intelligence, Toronto, ON, Canada}
	\affiliation{Chemical Physics Theory Group, Department of Chemistry, University of Toronto, Toronto, ON, Canada}
	\affiliation{Department of Computer Science, University of Toronto, Toronto, ON, Canada}
	\affiliation{Lebovic Fellow, Canadian Institute for Advanced Research, Toronto, ON, Canada}


\begin{abstract}
Numerous challenges in science and engineering can be framed as optimization tasks, including the maximization of reaction yields, the optimization of molecular and materials properties, and the fine-tuning of automated hardware protocols. Design of experiment and optimization algorithms are often adopted to solve these tasks efficiently. Increasingly, these experiment planning strategies are coupled with automated hardware to enable autonomous experimental platforms. The vast majority of the strategies used, however, do not consider robustness against the variability of experiment and process conditions. In fact, it is generally assumed that these parameters are exact and reproducible. Yet some experiments may have considerable noise associated with some of their conditions, and process parameters optimized under precise control may be applied in the future under variable operating conditions. In either scenario, the optimal solutions found might not be robust against input variability, affecting the reproducibility of results and returning suboptimal performance in practice. Here, we introduce \golem, an algorithm that is agnostic to the choice of experiment planning strategy and that enables robust experiment and process optimization. \golem identifies optimal solutions that are robust to input uncertainty, thus ensuring the reproducible performance of optimized experimental protocols and processes. It can be used to analyze the robustness of past experiments, or to guide experiment planning algorithms toward robust solutions on the fly. We assess the performance and domain of applicability of \golem through extensive benchmark studies and demonstrate its practical relevance by optimizing an analytical chemistry protocol under the presence of significant noise in its experimental conditions.
\end{abstract}
	\maketitle

\begin{bibunit}[unsrt]


\section{Introduction}

Optimization problems, in which one seeks a set of parameters that maximize or minimize an objective of interest, are ubiquitous across science and engineering. In chemistry, these parameters may be the experimental conditions that control the yield of the reaction, or those that determine the cost-efficiency of a manufacturing process (e.g., temperature, time, solvent, catalyst).\cite{Christensen:2020,Shields:2021} The design of molecules and materials with specific properties is also a multi-parameter, multi-objective optimization problem, with their chemical composition ultimately governing their properties.\cite{Nicolaou2013,Gomez-Bombarelli:2018,Sun:2019,Winter:2019,Yao:2021} These optimization tasks may, in principle, be performed autonomously. In fact, thanks to ever-growing automation, machine learning (ML)-driven experimentation has attracted considerable interest.\cite{Hase:2019,Gromski:2019,Stein:2019,Dimitrov:2019,Coley:2020a,Coley:2020b,FloresLeonar:2020} Self-driving laboratories are already accelerating the rate at which these problems can be solved by combining automated hardware with ML algorithms equipped with optimal decision-making capabilities.\cite{Nikolaev:2016,Maruyama:2017,Granda:2018,MacLeod:2020,Langner:2020,Grizou:2020,Tao:2021}

Recent efforts in algorithm development have focused on providing solutions to the requirements that arise from the practical application of self-driving laboratories. For instance, newly proposed algorithms include those with favorable computational scaling properties,\cite{Hase:2018} with the ability to optimize multiple objectives concurrently,\cite{Hase:2018_chimera} that are able to handle categorical variables (such as molecules) and integrate external information into the optimization process\cite{Hase:2020_gryffin}. One practical requirement of self-driving laboratories that has received little attention in this context is that of robustness against variability of experimental conditions and process parameters.

During an optimization campaign, it is typically assumed that the experimental conditions are known and exactly reproducible. However, the hardware (e.g., dispensers, thermostats) may impose limitations on the precision of the experimental procedure such that there is a stochastic error associated with some or all conditions. As a consequence, the optimal solution found might not be robust to perturbations of the inputs, affecting the reproducibility of the results and returning suboptimal performance in practice. Another scenario is when a process optimized under precise control is to be adopted in the future under looser operating conditions. For instance, in large-scale manufacturing, it might not be desirable (or possible) to impose tight operating ranges on the process parameters due to the cost of achieving high precision. This means that the tightly controlled input parameters used during optimization might not reflect the true, variable operating conditions that will be encountered in production.

In general, it is possible to identify two main types of input variability encountered in an experimental setting. The first is due to uncertainty in the experimental conditions that are controlled by the researchers, often referred to as the \textit{control factors}, corresponding to the examples discussed above. It can be caused by the imprecision of the instrumentation, which may reflect a fundamental limitation or a design choice, and could affect the present or future executions of the experimental protocol. A second type of input variability that can affect the performance of the optimization is due to experimental conditions that the researcher does not directly control. This may be, for instance, the temperature or the humidity of the room in which the experiments are being carried out. While it might not always be possible or desirable to control these conditions, they might be known and monitored such that their impact on the experimental outcome can in principle be accounted for.\cite{Pendleton:2019} The work presented here focuses on the first type of variability, related to control factors, although the approach presented may be in principle extended and applied to environmental factors too.

Here, we introduce \golem, a probabilistic approach that identifies optimal solutions that are robust to input uncertainty, thus ensuring the reproducible performance of optimized experiments and processes. \golem accounts for sources of uncertainty and may be applied to reweight the merits of previous experiments, or integrated into popular optimization algorithms to directly guide the optimization toward robust solutions. In fact, the approach is agnostic to the choice of experiment planning strategy and can be used in conjunction with both design of experiment and optimization algorithms. To achieve this, \golem explicitly models experimental uncertainty with suitable probability distributions that refine the merits of the collected measurements. This allows one to define an objective function that maximizes the average performance under variable conditions, while optionally also penalizing the expected variance of the results.

The article is organized as follows. First, we review some background information and previous work on robust optimization (section \ref{section:background}). Second, we introduce the core ideas behind the \golem algorithm (section \ref{section:golem}). We then present the analytical benchmark functions used to test \golem together with different optimization approaches (section \ref{section:surfaces}), as well as the results of these benchmark studies (section \ref{section:benchmarks}). Finally, we show how \golem may be used in practice, taking the calibration of a high-performance liquid chromatography (HPLC) protocol as an example application (section \ref{section:hplc}).


\section{Background and related work} 
\label{section:background}

Formally, an optimization task requires finding the set of conditions $\bm{x}$ (i.e., the \textit{parameters}, or \textit{control factors}) that yield the most desirable outcome for $f(\bm{x})$. If the most desirable outcome is the one that minimizes $f(\bm{x})$, then the solution of the optimization problem is

\begin{align}
\bm{x}^* = \argmin_{\bm{x} \in \mathcal{X}} f(\bm{x}),
\end{align}

where $\mathcal{X}$ is the domain of the optimization defining the range of experimental conditions that are feasible or that one is willing to consider. The objective function value $f(\bm{x})$ determines the \textit{merit} of a specific set of parameters $\bm{x}$. This merit may reflect the yield of a reaction, the cost-efficiency of a manufacturing process, or a property of interest for a molecule or material. Note that the objective function $f(\bm{x})$ is \textit{a priori} unknown, but can be probed via experiment. Only a finite number $K$ of samples $\mathcal{D}_{K} = \{\bm{x}, f(\bm{x})\}_{k=1}^{K}$ are typically collected during an optimization campaign, due to the cost and time of performing the experiments. A \textit{surrogate} model of $f(\bm{x})$ can be constructed based on $\mathcal{D}_{K}$. This model is typically a statistical or machine learning (ML) model that captures linear and non-linear relationships between the input conditions $\bm{x}$ and the objective function values $f(\bm{x})$.

An optimization campaign thus typically proceeds by iteratively testing sets of parameters $\bm{x}$, as defined via a design of experiment or as suggested by an experiment planning algorithm\cite{Hase:2020_olympus,Felton:2020,Rohr:2020}. Common design of experiment approaches rely on random or systematic searches of parameter combinations. Other experiment planning algorithms include sequential model-based approaches, such as Bayesian optimization\cite{Frazier:2018,Shahriari:2015}, and heuristic approaches like evolutionary and genetic algorithms\cite{Mitsuo:2008,McCall:2005,Srinivas:1994}. Experiment planning algorithms are now of particular interest in the context of self-driving laboratories for chemistry and materials science\cite{MacLeod:2020,Langner:2020,Hase:2018,Burger:2020,Gongora:2020}, which aim to autonomously and efficiently optimize the properties of molecules and materials.

\subsection{Robust optimization}

\begin{figure}[htb]
    \centering
    \includegraphics[width=1.0\columnwidth]{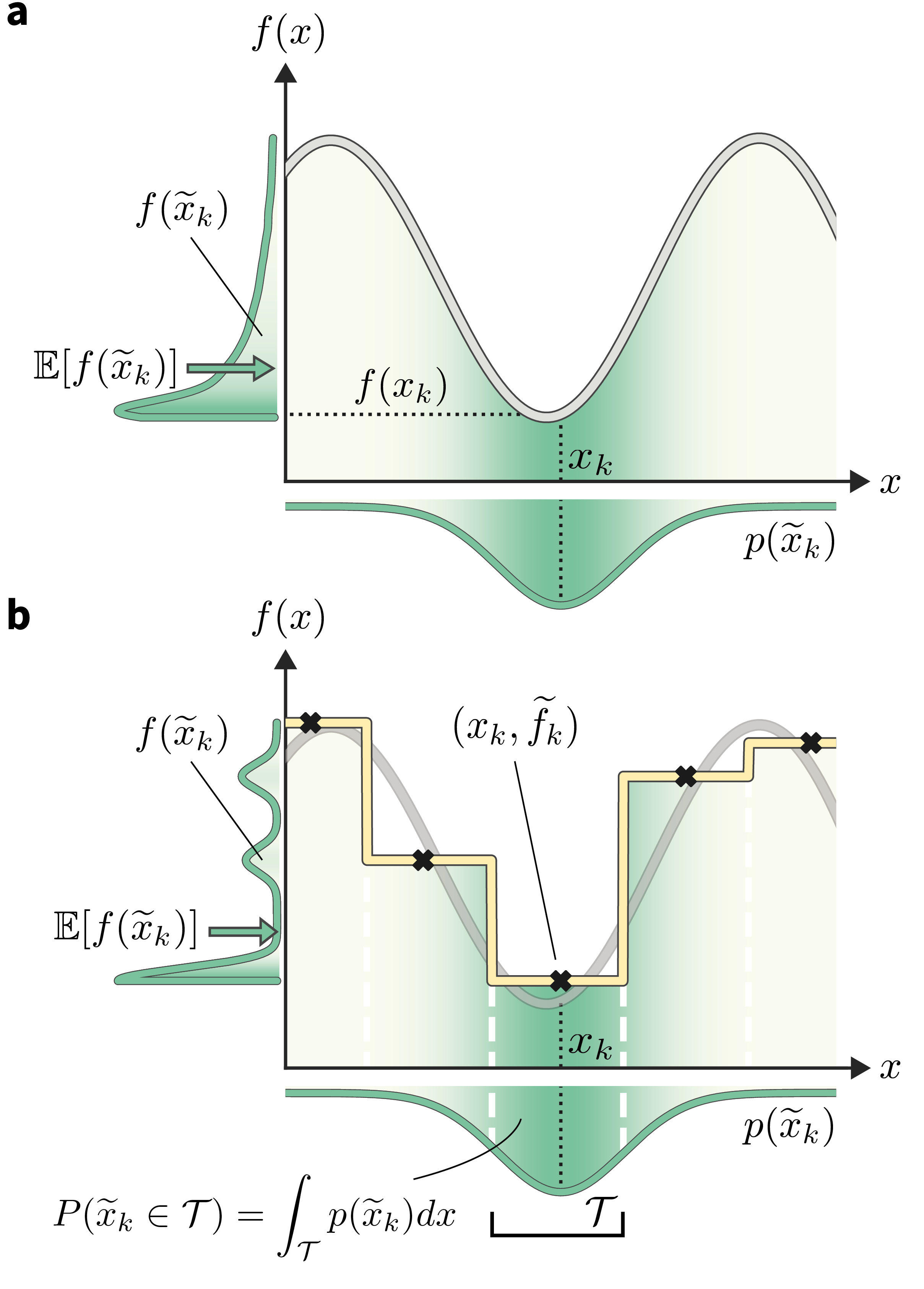}
    \caption{\golem's approach to estimating robustness. (\textbf{a}) Effect of uncertain inputs on objective function evaluations. The true objective function is shown as a gray line. The probability distribution $p(\widetilde{x}_k)$ of possible input value realizations for the targeted location $x_k$ is shown in green, below the $x$-axis. The distribution of output $f(\widetilde{x}_k)$ values caused by the input uncertainty are similarly shown next to the $y$-axis. The expectation of $f(\widetilde{x}_k)$ is indicated by a green arrow. (\textbf{b}) Schematic of \golem's core concept. The yellow line represents the surrogate function used to model the underlying objective function, shown in the background as a gray line. This surrogate is built with a regression tree, trained on five observations (black crosses). Note how the observations $\widetilde{f}_k$ are noisy, due to the uncertainty in the location of the input queries. In the noiseless query setting, and assuming no measurement error, the observations would lie exactly on the underlying objective function. Vertical white, dashed lines indicate how this model has partitioned the one-dimensional input space. Given a target location $x_k$, the probability that the realized input was obtained from partition $\mathcal{T}$ can be computed by integrating the probability density $p(\widetilde{x}_k)$ over $\mathcal{T}$, which is available analytically.}
    \label{fig:golem_scheme}
\end{figure}

The goal of \textit{robust} optimization is to identify solutions to an optimization problem that are robust to variation or sources of uncertainty in the conditions under which the experiments are or will be performed.\cite{Beyer:2007} Robustness may be sought for different reasons. For instance, the true location in parameter space of the query points being evaluated might be uncertain if experiments are carried out with imprecise instruments. In another scenario, a process might be developed in a tightly controlled experimental setting, however, it is expected that future execution of the same protocol will not. In such cases, a solution that is insensitive to the variability of the experimental conditions is desirable.

Several unique approaches have been developed for this purpose, originating with the robust design methodology of Taguchi, later refined by Box and others.\cite{Jones:2014,Beyer:2007} Currently, the most common approaches rely on either a \textit{deterministic} or \textit{probabilistic} treatment of input parameter uncertainty. Note that, by \textit{robust optimization}, and with chemistry applications in mind, we broadly refer to any approach aiming at solutions that mitigate the effects of the variability of experimental conditions. In the literature, the same term is sometimes used to specifically refer to what we are here referring to as \textit{deterministic} approaches.\cite{Bertsimas:2011,Beyer:2007} At the same time, the term \textit{stochastic optimization}\cite{Dantzig:1955,Powell:2019} is often used to refer to approaches that here we describe as \textit{probabilistic}. We also note that, while being separate fields, many similarities with robust control theory are present.\cite{Kemin:1996} The lack of a unified nomenclature is the result of robust optimization problems arising in different fields of science and engineering, from  operations research to robotics, finance, and medicine, each with their own sets of unique challenges. While a detailed review of all robust optimization approaches developed to date is out of the scope of this brief introductory section, we refer the interested reader to more comprehensive appraisals by Beyer\cite{Beyer:2007}, Bertsimas\cite{Bertsimas:2011}, and Powell\cite{Powell:2019}. In the interest of conciseness, we also do not discuss approaches based on fuzzy sets\cite{Chen2000,Inuiguchi2012} and those based on the minimization of risk measures\cite{Hong:2009,Cakmak2020}.

Deterministic approaches define robustness with respect to an uncertainty set.\cite{Bertsimas2009,Bogunovic:2018} Given the objective function $f(\bm{x})$, the robust counterpart $g(\bm{x})$ is defined as

\begin{align}
g(\bm{x}) \equiv \sup_{\bm{z} \in \mathcal{U}(\bm{x},\delta)} f(\bm{z}),
\end{align}

where $\mathcal{U}$ is an area of parameter space in the neighborhood of $\bm{x}$, the size of which is determined by $\delta$. $g(\bm{x})$ then takes the place of $f(\bm{x})$ in the optimization problem. This approach corresponds to optimizing for a worst-case scenario, since the robust merit is defined as the worst (i.e., maximum, in minimization tasks) value of $f(\bm{x})$ in the neighborhood of $\bm{x}$. Despite being computationally attractive, this approach is generally conservative and can result in robust solutions with poor average performance.\cite{Beyer:2007}

A different way to approach the problem is to treat input parameters probabilistically as random variables. Probability distributions for input parameters can be defined assuming knowledge about the uncertainty or expected variability of the experimental conditions.\cite{Beyer:2007} In this case, the objective function $f(\bm{x})$ becomes a random quantity itself, with its own (unknown) probability density (Figure \ref{fig:golem_scheme}a). The robust counterpart of $f(\bm{x})$ can then be defined as its expectation value,

\begin{align}
g(\bm{x}) \equiv \mathbb{E}[f(\widetilde{\bm{x}})] = \int f(\bm{x}) p({\bm{\widetilde{x}}}) d\bm{x}.
\label{eq:robust_prob}
\end{align}

Here, $\bm{\widetilde{x}} = \bm{{x}} + \bm{\delta}$, where $\bm{\delta}$ is a random variable with probability density $p(\bm{\delta})$, which represents the uncertainty of the input conditions at $\bm{x}$ (see section \ref{section:golem_derivation} for a different, but equivalent formulation). This definition ensures that the solution of the robust optimization problem is average-case optimal. For example, assume $f(\bm{x})$ is the yield of a reaction given the reaction conditions $\bm{x}$. However, we know the optimized protocol will be used multiple times in the future without carefully monitoring the experimental conditions. By optimizing $g(\bm{x})$ as defined above, instead of $f(\bm{x})$, and assuming that $p(\bm{\widetilde{x}})$ captures the variability of future experimental conditions correctly, one can identify a set of experimental conditions that returns the best possible yield on average across multiple repeated experiments.

Despite its attractiveness, the probabilistic approach to robust optimization presents computational challenges. In fact, the above expectation cannot be computed analytically for most combinations of $f(\bm{x})$ and $p(\bm{\widetilde{x}})$. One solution is to approximate $\mathbb{E}[f(\widetilde{\bm{x}})]$ by numerical integration, using quadrature or sampling approaches.\cite{Dellaportas1995,Pearce:2017,ToscanoPalmerin:2018} However, this strategy can become computationally expensive as the dimensionality of the problem increases and if $g(\bm{x})$ is to be computed for many samples. As an alternative numerical approach, it has been proposed to use a small number of carefully chosen points in $\bm{x}$ to cheaply approximate the integral.\cite{Nogueira:2016} Selecting optimal points for arbitrary probability distributions is not straightforward, however.\cite{Julier:2002}

In Bayesian optimization, it is common to use Gaussian process (GP) regression to build a surrogate model of the objective function. A few approaches have been proposed in this context to handle input uncertainty.\cite{Girard:2003,Damianou:2014} Most recently, Fr\"{o}hlich \textit{et al.}\cite{Frohlich:2020} have introduced an acquisition function for GP-based Bayesian optimization for the identification of robust optima. This formulation is analytically intractable and the authors propose two numerical approximation schemes. A similar approach was previously proposed by Beland and Nair\cite{Beland:2017}. However, in its traditional formulation, GP regression scales cubically with the number of samples collected. In practice, this means that optimizing $g(\bm{x})$ can become costly after collecting more than a few hundred samples. In addition, GPs do not inherently handle discrete or categorical variables\cite{Garrido-Merchan:2020} (e.g., type of catalyst), which are often encountered in practical chemical research. Finally, these approaches generally assume normally distributed input noise, as this tends to simplify the problem formulation. However, physical constraints on the experimental conditions may cause input uncertainty to deviate from this scenario, such that it would be preferable to be able to model any possible noise distribution.

In this work, we propose a simple, inexpensive, and flexible approach to probabilistic robust optimization. \golem enables the accurate modeling of experimental conditions and their variability for continuous, discrete, and categorical conditions, and for any (parametric) bounded or unbounded uncertainty distribution. By decoupling the estimation of the robust objective $g(\bm{x})$ from the details of the optimization algorithm, \golem can be used with any experiment planning strategy, from design of experiment, to evolutionary and Bayesian optimization approaches.


\section{Formulating Golem}
\label{section:golem}

\begin{figure}[tb]
    \centering
    \includegraphics[width=1.0\columnwidth]{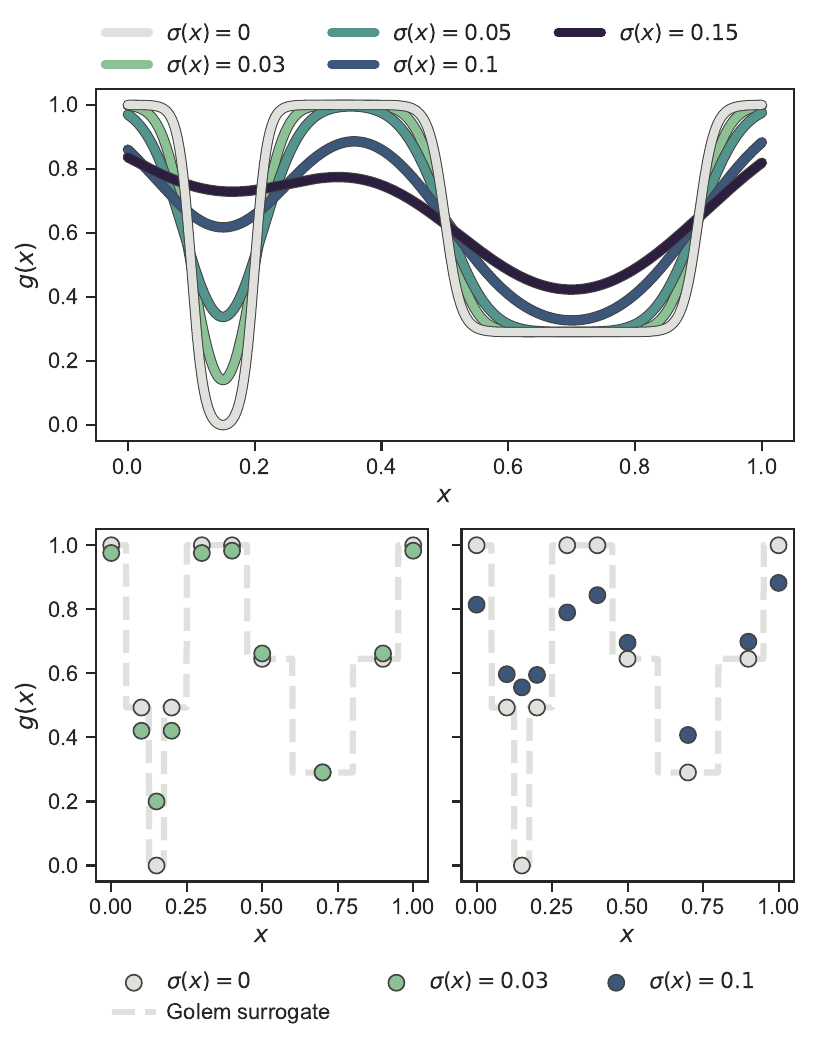}
    \caption{One-dimensional example illustrating the probabilistic approach to robustness and \golem's behavior. The top panel shows how $g(x)$, which is defined as $\mathbb{E}[f(x)] = \int f(x) p(\widetilde{x}) dx$, changes as the standard deviation of normally-distributed input noise $p(\widetilde{x})$ is increased. Note that the curve for $\sigma(x)=0$ corresponds to the original objective function. The panels at the bottom show the robust merits of a finite set of samples as estimated by \golem from the objective function values.}
    \label{fig:toy_example}
\end{figure}

Consider a robust optimization problem in which the goal is to find a set of input conditions $\bm{x} \in \mathcal{X}$ corresponding to the global minimum of the function $g : \mathcal{X} \to \mathbb{R}$,

\begin{align}
\bm{x}^* = \argmin_{\bm{x} \in \mathcal{X}} g(\bm{x}).
\label{eq:robust_opt_problem}
\end{align}

We refer to $g(\bm{x})$, as defined in Eq.~\ref{eq:robust_prob}, as the robust objective function, while noting that other integrated measures of robustness may also be defined. 

Assume a sequential optimization in which we query a set of conditions $\bm{x}_k$ at each iteration $k$. If the input conditions are noiseless, we can evaluate the objective function at $\bm{x}_k$ (denoted $f_k$). After $K$ iterations, we will have built a dataset $\mathcal{D}_{K} = \{ \bm{x}_k, f_k \}_{k=1}^K$. However, if the input conditions are noisy, the realized conditions will be $\widetilde{\bm{x}}_k = \bm{x}_k + \bm{\delta}$, where $\bm{\delta}$ is a random variable. As a consequence, we incur stochastic evaluations of the objective function, which we denote $\widetilde{f}_k$. This is illustrated in Figure \ref{fig:golem_scheme}a, where the Gaussian uncertainty in the inputs results in a broad distribution of possible output values. In this case, we will have built a dataset $\widetilde{\mathcal{D}} = \{ \bm{x}_k, \widetilde{f}_k \}_{k=1}^K$. Note that, while $\widetilde{\bm{x}}_k$ generally refers to a random variable, when considered as part of a dataset $\widetilde{\mathcal{D}}$ it may be interpreted as a specific sample of such variable. Hence, for added clarity, in Figure \ref{fig:golem_scheme} we refer to the distributions on the y-axis as $f(\widetilde{\bm{x}}_k)$, while we refer to function evaluations on specific input values as $\widetilde{f}_k$.

\subsection{General formalism}

The goal of \golem is to provide a simple and efficient means to estimate $g(\bm{x})$ from the available data, $\mathcal{D}$ or $\widetilde{\mathcal{D}}_{K}$. This would allow us to create a dataset $\mathcal{G}_{K} = \{ \bm{x}_k, g_k \} _{k=1}^K$ with robust merits, which can then be used to solve the robust optimization task in Eq.~\ref{eq:robust_opt_problem}. To do this, a surrogate model of the underlying objective function $f(\bm{x})$ is needed. This model should be able to capture complex, non-linear relationships. In addition, it should be computationally cheap to train and evaluate, and be scalable to high-data regimes. At the same time, we would like to flexibly model $p(\widetilde{\bm{x}})$, such that it can satisfy physical constraints and closely approximate the true experimental uncertainty. At the core of \golem is the simple observation that when approximating $f(\bm{x})$ with tree-based ML models, such as regression trees and random forest, estimates of $g(\bm{x})$ can be computed analytically as a finite series for any parametric probability density $p(\widetilde{\bm{x}})$. A detailed derivation can be found in section~\ref{section:golem_derivation}. 

An intuitive depiction of \golem is shown in Figure \ref{fig:golem_scheme}b. Tree-based models are piece-wise constant and rely on the rectangular partitioning of input space. Because of this discretization, $\mathbb{E}[f(\bm{x})]$ can be obtained as a constant contribution from each partition $\mathcal{T}$, weighted by the probability of $\bm{x}$ being within each partition, $P(\bm{x}_k \in \mathcal{T})$. Hence, an estimate of $g(\bm{x})$ can be efficiently obtained as a sum over all partitions (Eq.~\ref{eq:full_golem}).

Tree-based models such as regression trees and random forests have a number of advantages that make them well-suited for this task. First, they are non-linear ML models that have proved to be powerful function approximators. Second, they are fast to train and evaluate, adding little overhead to the computational protocols used. In the case of sequential optimization, the dataset $\mathcal{D}_{K}$ grows at each iteration $k$, such that the model needs to be continuously re-trained. Finally, they can naturally handle continuous, discrete, and categorical variables, so that uncertainty in all type of input conditions can be modeled. These reasons in addition to the fact that tree-based models allow for a closed-form solution to Eq.~\ref{eq:robust_prob} make \golem a simple yet effective approach for robust optimization. Note that while we decouple \golem's formulation from any specific optimization algorithm in this work, it is in principle possible to integrate this approach into tree-ensemble Bayesian optimization algorithms\cite{smac:2017,entmoot:2020}. This can be achieved via an acquisition function that is based on \golem's estimate of the robust objective, as well as its uncertainty, which can be estimated from the variance of $g(\bm{x})$ across trees.


Figure~\ref{fig:toy_example} shows a simple, one-dimensional example to provide intuition for \golem's behavior. In the top panel, the robust objective function is shown for different levels of normally-distributed input noise, parameterized by the standard deviation $\sigma(x)$ reported. Note that, when there is no uncertainty and $\sigma(x)=0$ (gray line), $p(x)$ is a delta function and one recovers the original objective function. As the uncertainty increases, the global minimum of the robust objective shifts from being the one at $x \approx 0.15$ to that at $x \approx 0.7$. In the two panels at the bottom, the same effect is shown under a realistic low-data scenario, in which only a few observations of the objective function are available (gray circles). Here, the dashed gray line represents the surrogate model used by \golem to estimate the robustness of each solution, given low (bottom left, green circles) and high (bottom right, blue circles) input noise. As in the top panel, which shows the continuous ground truth, here too the left-hand-side minimum is favored until the input noise is large enough such that the right-hand-side minimum provides better average-case performance.

\subsection{Multi-objective optimization}
When experimental noise is present, optimizing for the robust objective might not be the only goal. Often, large variance in the outcomes of an experimental procedure is undesirable, such that one might want to minimize it. For instance, in a chemical manufacturing scenario, one would like to ensure maximum overall output across multiple plants and batches. However, it would also be important that the amount of product manufactured in each batch does not vary considerably. Thus, the optimal set of manufacturing conditions should not only provide high yields on average, but also consistent ones. The problem can thus be framed as a multi-objective optimization in which we would like to maximize $\mathbb{E}[f(\bm{x})]$ while minimizing $\sigma[f(\bm{x})] = Var[f(\bm{x})]^{1/2}$. \golem can also estimate $\sigma[f(\bm{x})]$ (section~\ref{section:si_multiobj}), enabling such multi-objective optimizations. With $\mathbb{E}[f(\bm{x})]$ and $\sigma[f(\bm{x})]$ available, any scalarizing function may be used, including weighted sums and rank-based algorithms\cite{Hase:2018_chimera}.


\section{Benchmark surfaces and basic usage} 
\label{section:surfaces}

\begin{figure*}[htb]
    \centering
    \includegraphics[width=1.0\textwidth]{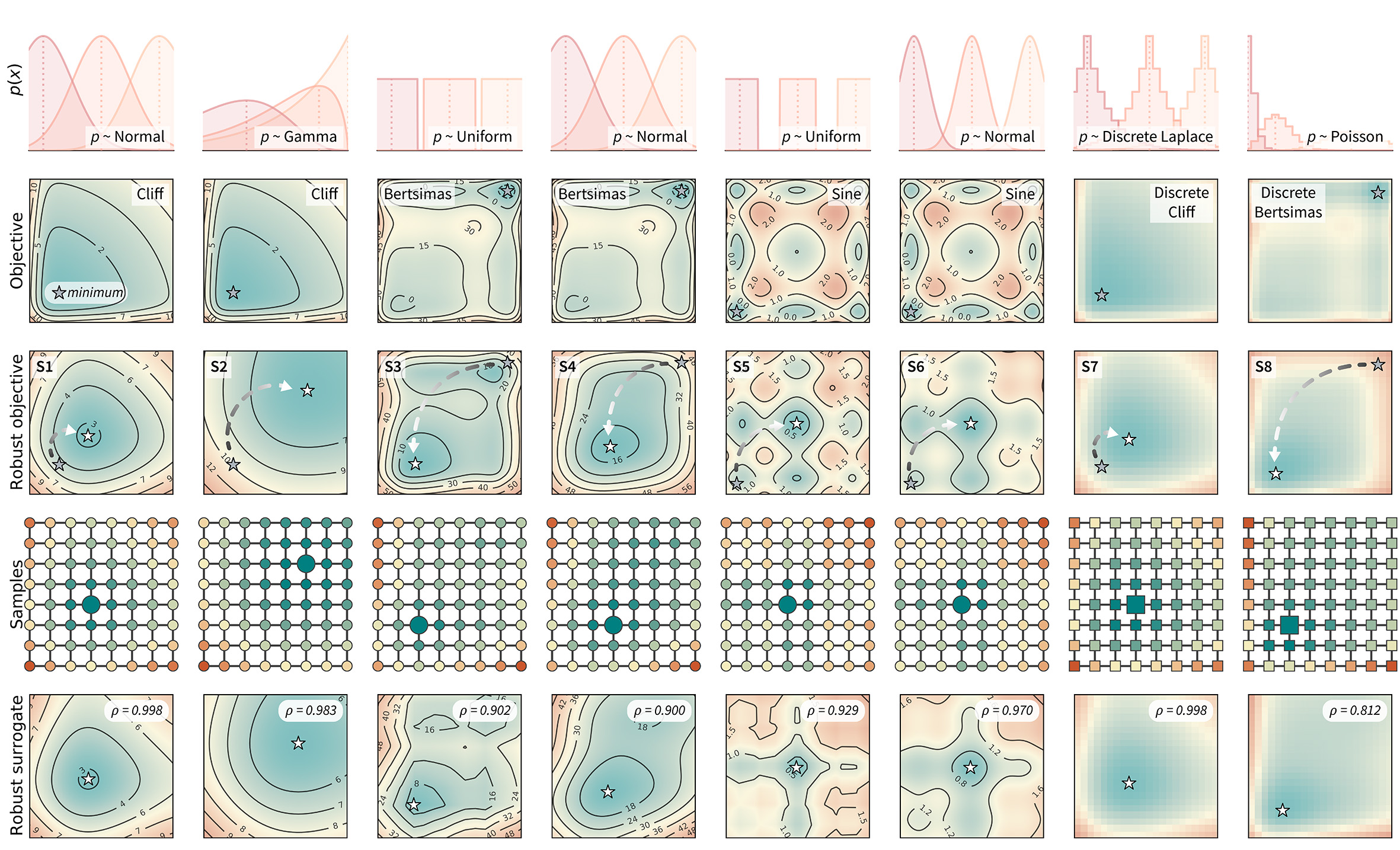}
    \caption{Benchmark functions used to test \golem and its performance. The first two rows show the type of uncertainty (in both input dimensions) assumed and the objective functions used in the synthetic benchmarks. The location of the global minimum is marked by a gray star on the two-dimensional surface of each objective function. The third row shows how the input uncertainty transforms each objective function into its robust counterpart. These surfaces (referred to as S1 to S8) represent the robust objectives, which are not directly observable, but that we would like to optimize. The global minimum of these functions are marked by white stars, with an arrow indicating the shift in the location of the global minimum between non-robust and robust objectives. The fourth row shows a set of $8\times 8$ samples that have been collected from these surfaces. Each sample is colored by its robust merit as estimated by \golem using only these $64$ samples. The larger marker (circle or square, for continuous and discrete surfaces, respectively) indicate the sample with best estimated robust merit. For all surfaces, \golem correctly estimates the most robust sample to be one in the vicinity of the true global minimum. The final row shows \golem's surrogate model of the robust objective, constructed from the grid of $64$ samples shown in row four. This surrogate model is highly correlated with the true underlying robust objective, as indicated by Spearman's correlation coefficient ($\rho$) reported at the top-right corner of each plot.}
    \label{fig:grids}
\end{figure*}

The performance of \golem, in conjunction with a number of popular optimization algorithms, was evaluated on a set of two-dimensional analytical benchmark functions. This allowed us to test the performance of the approach under different, hypothetical scenarios, test which optimization algorithms are most suited to be combined with \golem, and demonstrate the ways in which \golem may be deployed.

\subsection{Overview of the benchmark surfaces}

Figure~\ref{fig:grids} shows the benchmark functions that were used to evaluate \golem. These benchmarks were chosen to both challenge the algorithm and show its flexibility. We selected both continuous and discrete surfaces, and bounded and unbounded probability distributions to describe the input uncertainty. The objective functions considered are shown in the second row of Figure~\ref{fig:grids}. The \textit{Bertsimas} function is taken from the work of Bertsimas \textit{et al.}\cite{Bertsimas2009}, while \textit{Cliff} and \textit{Sine} are introduced in this work (section \ref{section:si_surfaces}). The first row of Figure~\ref{fig:grids} shows the uncertainty applied to these objective functions in both input dimensions. These uncertainties induce the robust objective functions shown in the third row. The location of the global minimum is shown for each objective and robust objective, highlighting how the location of the global minimum is affected by the variability of the inputs. The eight robust objectives in the third row of Figure~\ref{fig:grids} are labeled S1 to S8 and are the surfaces to be optimized. While we can only probe the objective functions in the second row, we use \golem to estimate their robust counterparts in the third row and locate their global minima.

These synthetic functions challenge \golem and the optimization algorithms in different ways. The rougher the surface and its robust counterpart, the more challenging it is expected to be to optimized. The smaller the difference in robust merit between the non-robust and robust minima (section \ref{section:si_surfaces}, Table~\ref{tab:surfaces}), the harder it is for \golem to resolve the location of the true robust minimum, as more accurate estimates of $g(\bm{x})$ are required. Finally, the steeper the objective function is outside the optimization domain, the less accurate \golem's estimate will be close to the optimization boundary, as samples are collected only within the optimization domain. 

S1--S6 evaluate performance on continuous spaces, while S7 and S8 on discrete ones. The function denoted \textit{Cliff} has a single minimum, which is shifted in the robust objectives S1 and S2. The \textit{Bertsimas} function has a global minimum indicated at the top-right corner of the surface, and a broader minimum at the bottom-left corner. The latter is the global minimum of the robust objective functions S3 and S4. The \textit{Sine} function is the most rugged and challenging, with nine minima (eight local and one global). S2 and S8 describe input uncertainty via distributions that do not allow values outside some of the bounds of the optimization domain. This is used to demonstrate \golem's flexibility and ability to satisfy physical constraints. For instance, if the uncertain input variable is dispensed volume, one should be able to assign zero probability to negative volumes.

\subsection{Reweighting previous results} 
One possible use of \golem is to reweight the merits of previously tested experimental conditions. Imagine, for instance, that we have accurately and precisely evaluated how temperature and catalyst concentration affect the yield of a reaction in the laboratory. To achieve this, we have performed $64$ experiments using a uniformly spaced $8 \times 8$ grid. Based on this data, we know which of the measured conditions provide the best yield. However, the same reaction will be used in other laboratories, or in larger-scale manufacturing, where these two variables will not be precisely controlled because, e.g., precise control is expensive or requires a complex experimental setup. Therefore, we would like to reweight the merit of each of the $64$ conditions previously tested, and identify which conditions are robust against variations in temperature and pressure. \golem allows one to easily compute these robust merits given the uncertainty in the input conditions. We tested \golem under this scenario and the results are shown in Figure~\ref{fig:grids}. In particular, the fourth row shows the grid of 64 samples taken from the objective function and reweighted with \golem. The color of each sample indicates their robust merit as estimated by \golem, with blue being more robust and red less robust. The largest marker indicates the sample estimated to have the best robust merit, which is in close proximity to the location of the true robust minimum for all surfaces considered.

Based on these $64$ samples, \golem can also build a surrogate model of the robust objective. This model is shown in the last row of Figure~\ref{fig:grids}. These estimates closely resemble the true robust surfaces in the third row. In fact, the Spearman's rank correlations ($\rho$) between \golem’s surrogates and the true robust objectives were $\geq 0.9$ for seven out of eight surfaces tested. For S8 only, while the estimated location of the global robust minimum was still correct, $\rho \approx 0.8$ due to boundary effects. In fact, while the robust objective depends also on the behavior of the objective function outside of the defined optimization domain, we sample the objective only within this domain. This lack of information causes the robustness estimates of points close to the boundaries to be less accurate than for those farther from them (Figure~\ref{si_fig:boundary_error}). Another consequence of this fact is that the robust surrogate does not exactly match the true robust objective also in the limit of infinite sampling within the optimization domain (section \ref{section:si_infsampling}).

To further clarify the above statement, by ``defined optimization domain'' we refer to a subset of the physically-meaningful domain that the researcher has decided to consider. Imagine, for instance, that we have a liquid dispenser which we will use to dispense a certain solvent volume. The smallest volume we can dispense is zero, while the largest might be the volume in the reservoir used (e.g., 1 L). These limits are physical bounds we cannot exceed. However, for practical purposes, we will likely consider a maximum volume much smaller than the physical limit (e.g., 5 mL). In this example, $0-5$ mL would constitute the defined optimization domain, while $0-1$ L are physical bounds on the domain. In the context of uncertain experimental conditions, it can thus be the case that a noisy dispenser might provide $5.1$ mL of liquid despite this exceeding the desired optimization boundary. The same cannot, however, be the case for the lower bound in this example, since a negative volume is physically impossible. As a consequence, while we allow an optimization algorithm to query the objective function only within the user-defined optimization domain, a noisy experimental protocol might result in the evaluation of the objective function outside of this domain.

\golem allows to take physical bounds into account by modeling input uncertainty with bounded probability distributions. Yet, it cannot prevent boundary effects that are the consequence of the unknown behaviour of the objective function outside of the defined optimization domain. This issue, unfortunately, cannot be resolved in a general fashion, as it would require a data-driven model able to extrapolate arbitrarily far from the data used for training. A practical solution may be to consider a ``data collection domain'' as a superset of the optimization domain, which is used for collecting data at the boundaries but which the optimization solution is not selected from. In the examples in Figure~\ref{fig:grids} (row 4), this would mean using the datapoints on the perimeter of the two-dimensional grid only for estimating the robustness of the internal points more accurately. We conclude by reiterating how, notwithstanding this inescapable boundary effect, as shown in Figure~\ref{fig:grids} there is a high correlation between \golem's estimates and the true robustness values. 

\section{Optimization Benchmarks} 
\label{section:benchmarks}

\begin{figure*}[htb]
    \centering
    \includegraphics[width=1.0\textwidth]{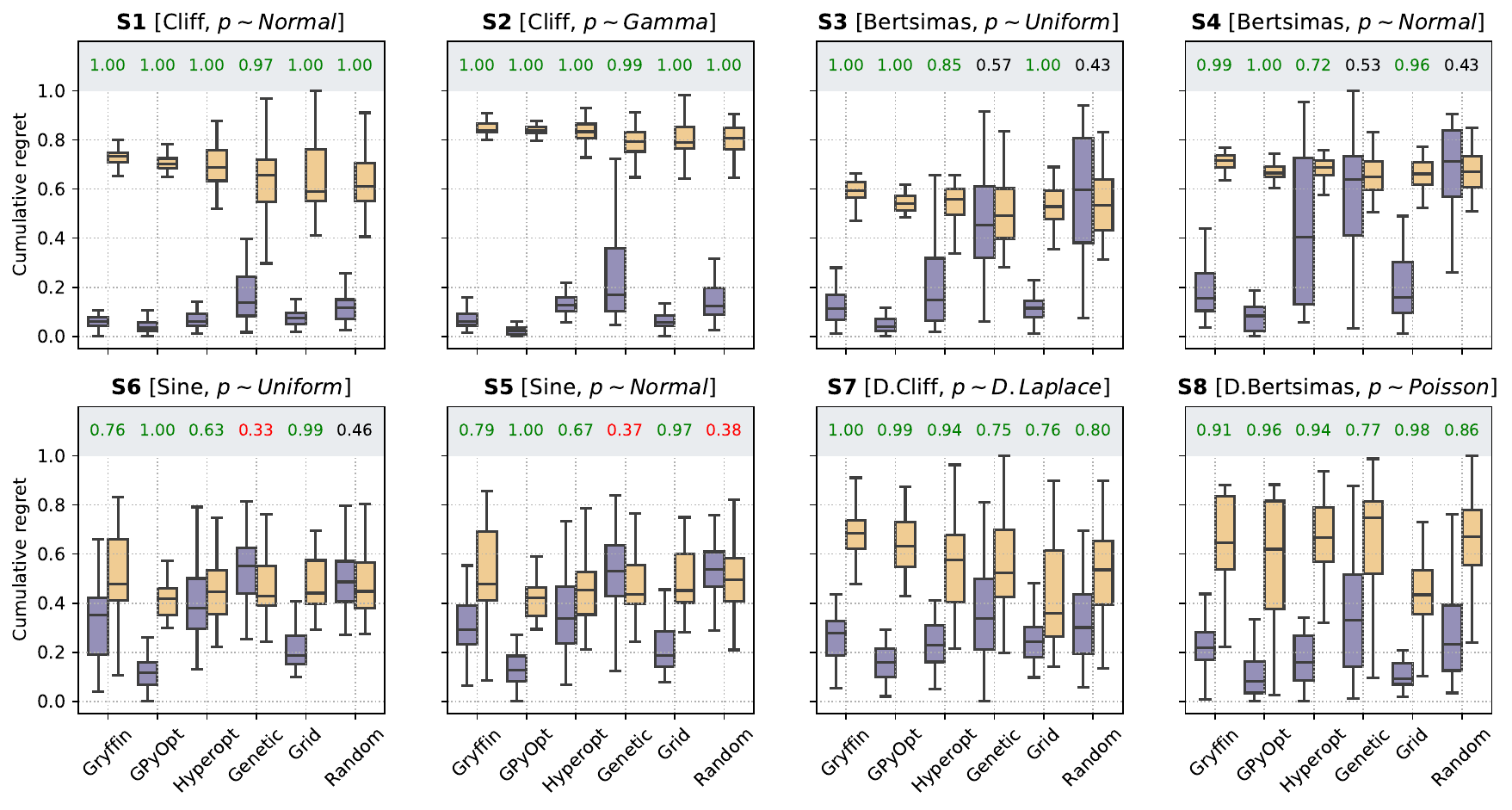}
    \caption{Robust optimization performance of multiple algorithms, with and without \golem, in benchmarks where queries were noiseless. Box plots show the distributions of cumulative regrets obtained across 50 optimization repeats with and without \golem, in purple and yellow, respectively. The boxes show the first, second, and third quartiles of the data, with whiskers extending up to 1.5 times the interquartile range. At the top of each plot, we report the probability that the use of \golem improved upon the performance of each algorithm. Probabilities are in green if the performance with \golem was significantly better (considering a 0.05 significance level) than without, and in red if it was significantly worse, as computed by bootstrap.}
    \label{fig:noiseless_benchmarks}
\end{figure*}

With increasing levels of automation and interest in self-driving laboratories, sequential approaches that make use of all data collected to select the next, most informative experiment are becoming the methods of choice for early prototypes of autonomous science. In this case, rather than re-evaluating previously performed experiments, one would like to steer the optimization towards robust solutions during the experimental campaign. \golem allows for this in combination with popular optimization approaches, by mapping objective function evaluations onto an estimate of their robust merits at each iteration of the optimization procedure. We evaluated the ability of six different optimization approaches to identify robust solutions when used with \golem and without. The algorithms tested include three Bayesian optimization approaches (\textit{Gryffin}\cite{Hase:2018,Hase:2020_gryffin}, \textit{GPyOpt}\cite{Gpyopt:2016}, \textit{Hyperopt}\cite{Bergstra:2013}), a genetic algorithm (\textit{Genetic})\cite{deap}, a random sampler (\textit{Random}), and a systematic search (\textit{Grid}). \textit{Gryffin}, \textit{GPyOpt}, and \textit{Hyperopt} use all previously collected data to decide which set of parameters to query next, \textit{Genetic} uses part of the collected data, while \textit{Random} and \textit{Grid} are totally agnostic to previous measurements.

In these benchmarks, we allowed the algorithms to collect $196$ samples for continuous surfaces and $64$ for the discrete ones. We repeated each optimization $50$ times to collect statistics. For \textit{Grid}, we created a set of $14 \times 14$ uniformly-spaced samples ($8 \times 8$ for the discrete surfaces) and then selected them at random at each iteration. For all algorithms tested, we performed the optimization with and without \golem. Algorithm performance in the absence of \golem constitutes a na\"ive baseline. Optimization performance in quantified using normalized cumulative robust regret, defined in \ref{section:si_regret}. This regret is a relative measure of how fast each algorithm identifies increasingly robust solutions, allowing the comparison of algorithm performance with respect to a specific benchmark function.

\subsection{Noiseless queries with uncertainty in future experiments}

Here, we tested \golem under a scenario where queries during the optimization are deterministic, i.e., noiseless. It is assumed that uncertainty in the inputs will arise only in future experiments. This scenario generally applies to the development of experimental protocols that are expected to be repeated under loose control of experimental conditions.

The results of the optimization benchmarks under this scenario are summarized in Figure \ref{fig:noiseless_benchmarks}, which shows the distributions of cumulative regrets for all algorithms considered, with and without \golem, across the eight benchmark surfaces. For each algorithm, Figure \ref{fig:noiseless_benchmarks} also quantifies the probability that the use of \golem resulted in better performance in the identification robust solutions. Overall, these results showed that \golem allowed the optimization algorithms to identify solutions that were more robust than those identified without \golem. 

A few additional trends can be extracted from Figure \ref{fig:noiseless_benchmarks}. The Bayesian optimization algorithms (\textit{Gryffin}, \textit{GPyOpt}, \textit{Hyperopt}) and systematic searches (\textit{Grid}) seemed to benefit more from the use of \golem than genetic algorithms (\textit{Genetic}) and random searches (\textit{Random}). In fact, the former approaches benefited from \golem across all benchmark functions, while the latter did so only for half the benchmarks. The better performance of \textit{Grid} as compared to \textit{Random}, in particular, may appear surprising. We found that the main determinant of this difference is the fact that \textit{Grid} samples the boundaries of the optimization domain, while \textit{Random} is unlikely to do so. By forcing random to sample the optimization boundaries, we recovered performances comparable to \textit{Grid} (section \ref{section:si_lowdisc}). We also hypothesized that uniformity of sampling might be beneficial to \golem, given that the accuracy of the robustness estimate depends on how well the objective function is modeled in the vicinity of the input location considered. We indeed found that low-discrepancy sequences provided, in some cases, slightly better performance than random sampling. However, this effect was minor compared to that of forcing the sampling of the optimization domain boundaries (section \ref{section:si_lowdisc}). 

\textit{Genetic} likely suffered from the same pathology, given it is initialized with random samples. Thus, in this context, initialization with a grid may be more appropriate. Genetic algorithms are also likely to suffer from a second effect. Given that we can only \textit{estimate} the robust objective, \golem induces a history-dependent objective function. Contrary to Bayesian optimization approaches, genetic algorithms consider only a subset of the data collected during optimization, as they discard solutions with bad fitness. Given that the robustness estimates change during the course of the optimization, these algorithms may drop promising solutions early in the search, which are then not recovered in the latter stages when \golem would have more accurately estimated their robustness. The use of more complex genetic algorithm formulations, exploring a more diverse set of possible solutions\cite{Nigam:2020}, could improve this scenario and is a possibility left for future work.

\subsection{Noisy queries with uncertainty in current experiments}
\label{section:benchmarks_noisy}

\begin{figure*}[htb]
    \centering
    \includegraphics[width=1.0\textwidth]{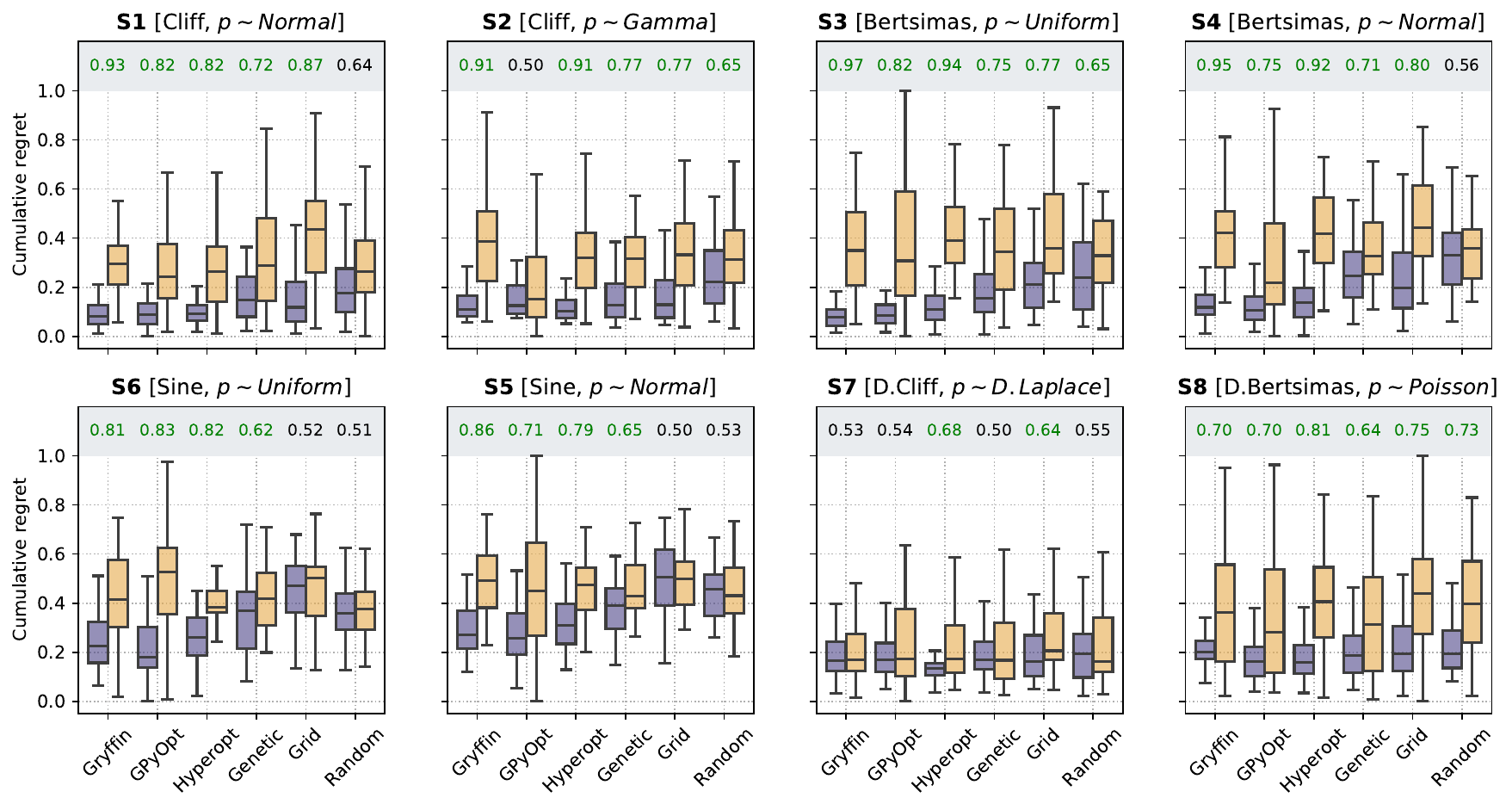}
    \caption{Robust optimization performance of multiple algorithms, with and without \golem, in benchmarks where queries were noisy. Box plots show the distributions of cumulative regrets obtained across 50 optimization repeats with and without \golem, in purple and yellow, respectively. The boxes show the first, second, and third quartiles of the data, with whiskers extending up to 1.5 times the interquartile range. At the top of each plot, we report the probability that the use of \golem improved the performance of each algorithm. Probabilities are in green if the performance with \golem was significantly better (considering a 0.05 significance level) than without, and in red if it was significantly worse, as computed by bootstrap.}
    \label{fig:noisy_benchmarks}
\end{figure*}

In a second scenario, queries during the optimization are stochastic, i.e., noisy, due the presence of substantial uncertainty in the current experimental conditions. This case applies to any optimization campaign in which it is not possible to precisely control the experimental conditions. However, we assume one can model the uncertainty $p(\widetilde{\bm{x}})$, at least approximately. For instance, this uncertainty might be caused by some apparatus (e.g., a solid dispenser) that is imprecise, but can be calibrated and the resulting uncertainty quantified. The optimization performances of the algorithms considered, with and without \golem, are shown in Figure \ref{fig:noisy_benchmarks}. Note that, to model the robust objective exactly, $p(\widetilde{\bm{x}})$ should also be known exactly. While this is not a necessary assumption of the approach, the accuracy of \golem's estimates is proportional to the accuracy of the $p(\widetilde{\bm{x}})$ estimates. As the $p(\widetilde{\bm{x}})$ estimate provided to \golem deviates from its true values, \golem under- or over-estimate the robustness of the optimal solution, depending on whether the input uncertainty is  under- or over-estimated. We will illustrate this point in more detail in section~\ref{section:hplc_retro}.

Generally speaking, this is a more challenging scenario than when queries are noiseless. As a consequence of the noisy experimental conditions, the dataset collected does not correctly match the realized control factors $\bm{x}$ with their associated merit $f(\bm{x})$. Hence, the surrogate model is likely to be a worse approximation of the underlying objective function than when queries are noiseless. While the development of ML models capable of recovering the objective function $f(\bm{x})$ based on noisy queries $\widetilde{\bm{x}}$ is outside the scope of this work, such models may enable even more accurate estimates of robustness with \golem. We are not aware of approaches capable of performing such an operation, but it is a promising direction for future research. In fact, being able to recover the (noiseless) objective function from a small number of noisy samples $\widetilde{f}$ would be beneficial not only for robustness estimation, but for the interpretation of experimental data more broadly.

Because of the above-mentioned challenge in the contruction of an accurate surrogate model, in some cases, the advantage of using \golem might not seem as stark as in the noiseless setting. This effect may be seen in surfaces S1 and S2, where the separation of the cumulative regret distributions is larger in Figure \ref{fig:noiseless_benchmarks} than it is in Figure \ref{fig:noisy_benchmarks}. Nonetheless, across all benchmark functions and algorithms considered, the use of \golem was beneficial in the identification of robust solutions in the majority of cases, and never detrimental, as shown by Figure \ref{fig:noisy_benchmarks}. In fact, \golem appears to be able to recover significant correlations with the true robust objectives $g(\bm{x})$ even when correlation with the objective functions $f(\bm{x})$ is lost due to noise the queried locations (Figure~\ref{si_fig:corr_fx_vs_gx}).

\begin{figure}[htb]
    \centering
    \includegraphics[width=0.9\columnwidth]{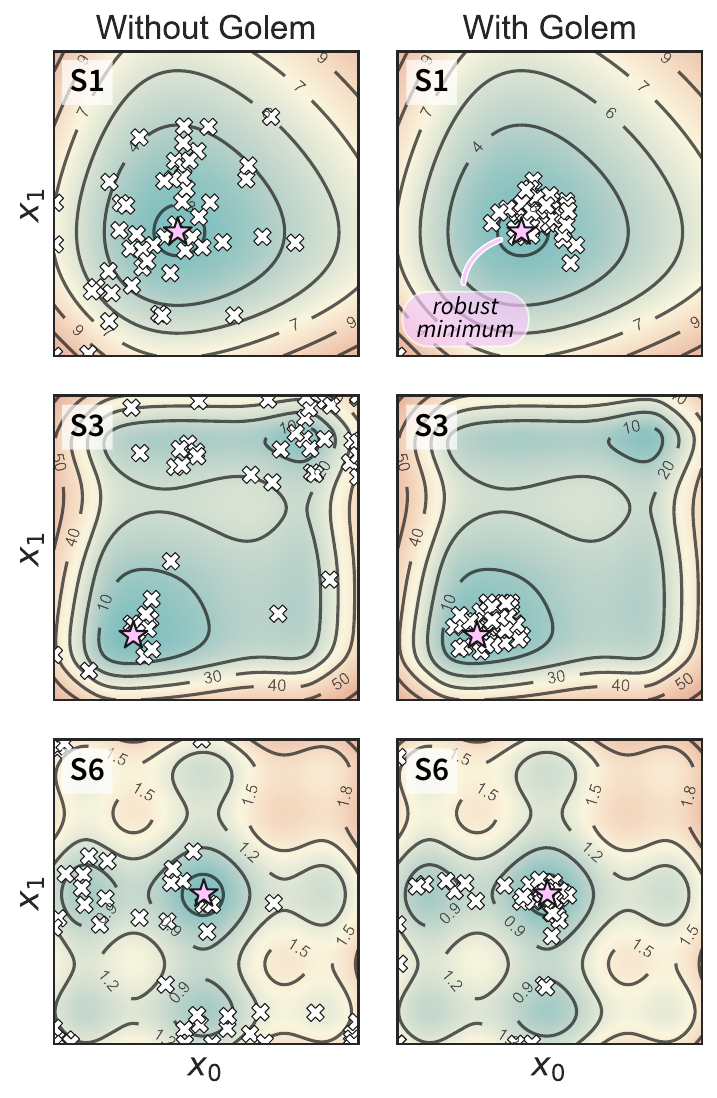}
    \caption{Location of the optimal input parameters identified with and without \golem. The results shown were obtained with \textit{GPyOpt} as the optimization algorithm. A pink star indicates the location of the true robust minimum. White crosses (one per optimization repeat, for a total of 50) indicate the locations of the optimal conditions identified by the algorithm without (on the left) and with (on the right) \golem.}
    \label{fig:xmin_examples}
\end{figure}

Optimization with noisy conditions is significantly more challenging than traditional optimization tasks with no input uncertainty. However, the synthetic benchmarks carried out suggest that \golem is able to efficiently guide optimization campaigns towards robust solutions. For example, Figure \ref{fig:xmin_examples} shows the location of the best input conditions as identified by \textit{GPyOpt} with and without \golem. Given the significant noise present, without \golem, the optima identified by different repeated experiments are scattered far away from the robust minimum. When \golem is used, the optima identified are considerably more clustered around the robust minimum.

\subsection{Effect of forest size and higher input dimensions}
\label{section:benchmarks_highdim}
All results shown thus far were obtained using a single regression tree as \golem's surrogate model. However, \golem can also use tree-ensemble approaches, such as random forest\cite{Breiman:2001} and extremely randomized trees\cite{Geurts:2006}. We thus repeated the synthetic benchmarks discussed above using these two ML models, with forest sizes of $10$, $20$, and $50$ (section \ref{section:si_forests}). Overall, for these two-dimensional benchmarks we did not observe significant improvements when using larger forest sizes. For the benchmarks in the noiseless setting, regression trees appeared to provide slightly better performance against the \textit{Bertsimas} functions (Figure \ref{si_fig:tree_dependance1}). The lack of regularization may have provided a small advantage in this case, where \golem is trying to resolve subtle differences between competing minima. Yet, a single regression tree performed as well as ensembles. For the benchmarks in the noisy setting, random forest and extremely randomized trees performed slightly better overall (Figure \ref{si_fig:tree_dependance2}). However, larger forests did not appear to provide considerable advantage over smaller ones, suggesting that for these low-dimensional problems, small forests or even single trees can generally be sufficient.

To study the performance of different tree-ensemble approaches also on higher-dimensional search spaces, we conducted experiments, similar to the ones described above, on three-, four-, five, and six-dimensional versions of benchmark surface S1. In these tests, we consider two dimensions to be uncertain, while the additional dimensions are noiseless. Here, too, we studied the effect of forest type and size on the results, but we focused on the Bayesian optimization algorithms. In this case, we observed better performance of \golem when using random forest or extremely randomized trees as the surrogate model. In the noiseless setting, extremely randomized trees returned slightly better performance than random forest, in particular for \textit{GPyOpt} and \textit{Hyperopt} (Figure \ref{si_fig:highdim_noiseless}). The correlation of optimization performance with forest size was weaker. Yet, for each combination of optimization algorithms and benchmark surface, the best overall performance was typically achieved with larger forest sizes of $20$ or $50$ trees. While less marked, similar trends were observed for the same tests in the noisy setting (Figure \ref{si_fig:highdim_noisy}). In this scenario, random forest returned slightly better performance than extremely randomized trees for \textit{Hyperopt}. Overall, surrogate models based on random forest or extremely randomized trees appear to provide better performance across different scenarios.

We then investigated \golem's performance across varying search space dimensionality and number of uncertain conditions. To do this, we conducted experiments on three-, four-, five, and six-dimensional versions of benchmark surface S1, with one to six uncertain inputs. These tests showed that \golem was still able to guide the optimizations towards better robust solutions. In the noiseless setting, the performance of \textit{GPyOpt} and \textit{Hyperopt} was significantly better with \golem for all dimensions and number of uncertain variables tested (Figure \ref{si_fig:highdim_highunc_noiseless}). The performance of \textit{Gryffin} was significantly improved by \golem in roughly half of the cases. Overall, given a certain search space dimensionality, the positive effect of \golem became more marked with a higher number of uncertain inputs. This observation does not imply that the optimization task is easier with more uncertain inputs (it is in fact more challenging), but that the use of \golem provides a more significant advantage in such scenarios. On the contrary, given a specific number of uncertain inputs, the effect of \golem was less evident with increasing number of input dimensions. Indeed, additional input dimensions make it more challenging for \golem to resolve whether the observed variability in the objective function evaluations is due to the uncertain variables or the expected behavior of the objective function along the additional dimensions. Similar overall results were observed in the noisy input setting (Figure \ref{si_fig:highdim_highunc_noisy}). However, statistically significant improvements were found in a smaller fraction of cases. Here, we did not observe a significant benefit in using \golem when having a small ($1-2$) number of uncertain inputs, but this became more evident with a larger ($3-6$) number of uncertain inputs. In fact, the same trends with respect to the dimensionality of the search space and the number of uncertain inputs were observed also in the noisy query setting. One important observation is that \golem was almost never (one out of $108$ tests) found to be detrimental to optimization performance, suggesting that there is very little risk in using the approach when input uncertainty is present, as in the worst-case scenario \golem would simply leave the performance of the optimization algorithm used unaltered.

Overall, these results suggest that \golem is also effective on higher-dimensional surfaces. In addition, it was found that the use of surrogate models based on forests can, in some cases, provide a better optimization performance. Given the limited computational cost of \golem, we thus generally recommend the use of an ensemble tree method as the surrogate model. Forest sizes of $20$ to $50$ trees were found to be effective. Yet, given that larger ensembles will not negatively affect the estimator performance, and that the runtime scales linearly with the number of trees, larger forests may be used as well.


\section{Chemistry Applications}
\label{section:hplc}

In this section, we provide an example application of \golem in chemistry. Specifically, we consider the calibration of an HPLC protocol, in which six controllable parameters (Figure \ref{fig:hplc1}a, section \ref{section:si_hplc}) can be varied to maximize the peak area, i.e., the amount of drawn sample reaching the detector.\cite{Roch:2020,Hase:2020_olympus} Imagine we ran $1386$ experiments in which we tested combinations of these six parameters at random. The experiment with the largest peak area provides the best set of parameters found. The parameter values corresponding to this optimum are highlighted in Figure \ref{fig:hplc1}b by a gray triangle pointing towards the abscissa. With the collected data, we can build a surrogate model of the response surface. The one shown as a gray line in Figure \ref{fig:hplc1}b was built with $200$ extremely randomized trees\cite{Geurts:2006}. Figure \ref{fig:hplc1}b shows the predicted peak area when varying each of the six controllable parameters independently around the optimum identified.

\subsection{Analysis of prior experimental results}
\label{section:hplc_retro}

\begin{figure}[tb]
    \centering
    \includegraphics[width=0.84\columnwidth]{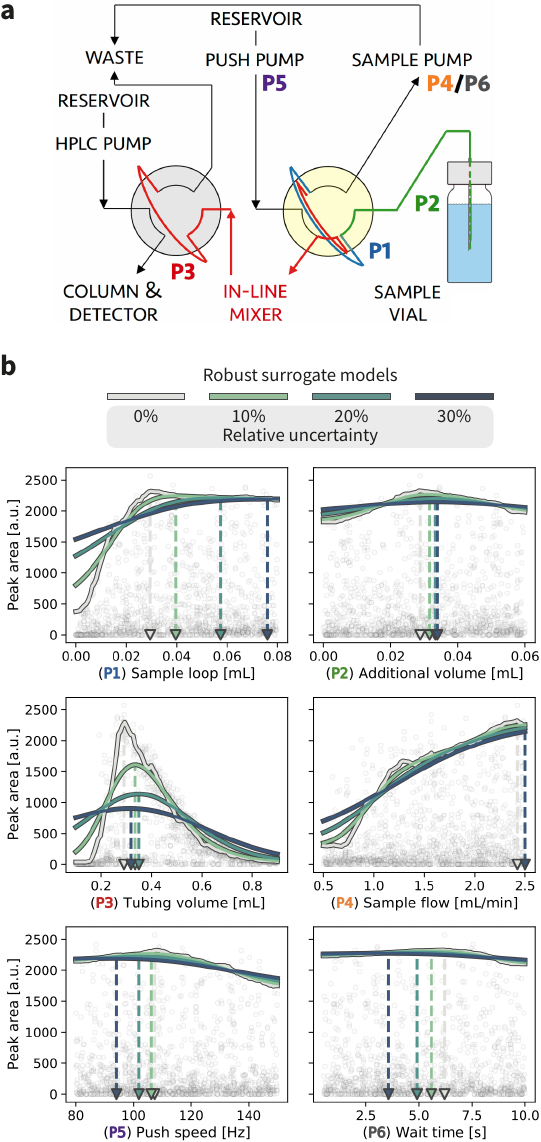}
    \caption{Analysis of the robustness of an HPLC calibration protocol. (a) Flow path for the HPLC sampling sequence performed by a robotic platform. The six parameters (P1-P6) are color coded. The yellow shade highlights the arm valve, and the gray shade the HPLC valve. (b) \golem analysis of the effect of input noise on expected protocol performance. A surrogate model of the response surface is shown in gray. Uncertainties were modeled with truncated normal distributions with standard deviations of 10\%, 20\%, 30\% of each parameter's range. The corresponding robust surrogate models are shown in light green, dark green, and blue. Triangular markers and dashed lines indicate the location of the optima for each parameter under different levels of noise.}
    \label{fig:hplc1}
\end{figure}

\golem allows us to speculate how the expected performance of this HPLC protocol would be affected by varying levels of noise in the controllable parameters. We modeled input noise via truncated normal distributions that do not support values below zero. This choice satisfies the physical constraints of the experiment, given that negative volumes, flows, and times are not possible. We considered relative uncertainties corresponding to a standard deviation of $10$\%, $20$\%, and $30$\% of the allowed range for each input parameter. The protocol performance is most affected by uncertainty in the tubing volume (variable P3, Figure~\ref{fig:hplc1}b). A relative noise of 10\% would result in an average peak area of around $1500$ a.u., a significant drop from the maximum observed at over $2000$. It follows that to achieve consistent high performance with this protocol, efforts should be spent in improving the precision of this variable. 

While the protocol performance (i.e., expected peak area) is least robust against uncertainty in P3, the location of the optimum setting for P3 is not particularly affected. Presence of noise in the sample loop (variable P1) has a larger effect on the location of its optimal settings. In fact, noise in P1 requires larger volumes to be drawn into the sample loop to be able to achieve average optimal responses. The optimal parameter settings for the push speed (P5) and wait time (P6) are also affected by the presence of noise. However, the protocol performance is fairly insensitive to changes in these variables, with expected peak areas of around $2000$ a.u. for any of their values within the range studied. 

Figure~\ref{fig:hplc1} also illustrates the effect of under- or over-estimating experimental condition uncertainty on \golem's robustness estimates. Imagine that the true uncertainty in variable P3 is $20$\%. This may be the true uncertainty encountered in the future deployment of the protocol, or it may be the uncertainty encountered while trying to optimize it. If we assume, incorrectly, the uncertainty to be $10$\%, \golem will predict the protocol to return, on average, an area of $\sim 1500$ a.u., while we will find that the true average performance of the protocol provides an area slightly above $1000$ a.u. That is, \golem will overestimate the robustness of the protocol. On the other hand, if we assumed the uncertainty to be $30$\%, we would underestimate the robustness of the protocol, as we would expect an average area below $1000$ a.u. In the case of variable P3, however, the location of the optimum is only slightly affected by uncertainty, such that despite the incorrect prediction, \golem would still accurately identify the location of the global optimum. That is, a tubing volume of $\sim 0.3$ mL provides the best average outcome whether the true uncertainty is $10$\%, $20$\%, or $30$\%. In fact, while ignoring uncertainty altogether (i.e. assuming $0$\% uncertainty) would result in the largest overestimate of robustness, it would still have minimal impact in practice given that the prediction of the optimum location would still be accurate. This is not the case if we considered P1. If we again assume that the true uncertainty in this variable is $20$\%, providing \golem with an uncertainty model with $10$\% standard deviation would result in a protocol using a sample loop volume of $\sim 0.04$ mL, while the optimal one should be $\sim 0.06$ mL. Providing \golem with a $30$\% uncertainty instead would result in an underestimate of the protocol robustness and an unnecessarily conservative choice of $\sim 0.08$ mL as the sample loop volume.

In summary, as anticipated in section~\ref{section:benchmarks_noisy}, while an approximate estimate of $p(\widetilde{\bm{x}})$ does not prevent the use of \golem, it can affect the quality of its predictions. When uncertainty is underestimated, the optimization solutions identified by \golem will tend to be less robust than expected. On the contrary, when uncertainty is overestimated, \golem's solutions will tend to be overly conservative (i.e., \golem will favor plateaus in the objective function despite more peaked optima would provide better average performance). The errors in \golem's estimates will be proportional to the error in the estimates of the input uncertainty provided to it, but the magnitude of these errors is difficult to predict as it depends on the objective function, which is unknown and application-specific. Note that, ignoring input uncertainty corresponds to assuming $p(\widetilde{\bm{x}})$ is a delta function in \golem. This choice, whether implicitly or explicitly made, results in the largest possible overestimate of robustness when uncertainty is in fact present. The associated error in the expected robustness is likely to be small when the true uncertainty is small, but may be large otherwise.

It is important to note that, above, we analyzed only one-dimensional slices of the six-dimensional parameter space. Given interactions between these parameters, noise in one parameter can affect the optimal setting of a different one (section \ref{section:si_hplc_interac}). \golem can identify these effects by studying its multi-dimensional robust surrogate model. Furthermore, for simplicity, here we considered noise in each of the six controllable parameters one at a time. It is nevertheless possible to consider concurrent noise in as many parameters as desired.

This example shows how \golem may be used to analyze prior experimental results and study the effect of input noise on protocol performance and the optimal setting of its controllable parameters.

\subsection{Optimization of a noisy HPLC protocol}
\label{section:hplc_opt}

\begin{figure*}[tb]
    \centering
    \includegraphics[width=1.0\textwidth]{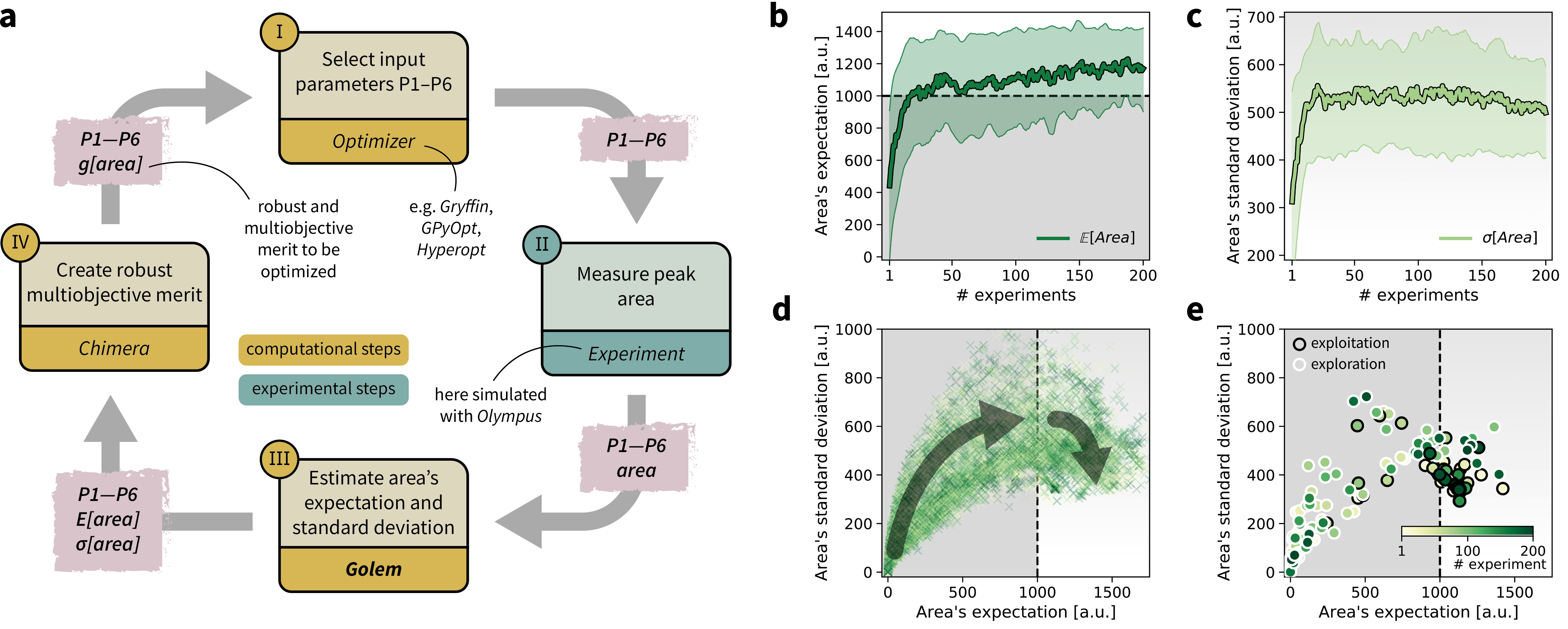}
    \caption{Setup and results for the optimization of an HPLC protocol under noisy experimental conditions. 
    (a) Procedure and algorithms used for the robust optimization of the HPLC protocol. First, the optimization algorithm selects the conditions of the next experiment to be performed. Second, the HPLC experiment is carried out and the associated peak's area recorded. Note that, in this example, P1 and P3 are noisy such that their values realized in the experiment do not correspond to those requested by the optimizer. Third, \golem is used to estimate the expected peak's area, $\mathbb{E}[Area]$, as well as its variability $\sigma[Area]$, based on a model of input noise for P1 and P3. Finally, the \textit{Chimera} scalarizing function is used to combine these two objectives into a single figure of merit to be optimized.
    (b-e) Results of 50 optimization repeats performed with \textit{Gryffin}. Equivalent results obtained with \textit{GPyOpt} and \textit{Hyperopt} are shown in Figure \ref{si_fig:hplc_opt_gpyopt_hopt}.
    (b) Optimization trace for the primary objective, i.e. the maximization of $\mathbb{E}[Area]$ above $1000$ a.u. The average and standard deviation across $50$ optimization repeats are shown.
    (c) Optimization trace for the secondary objective, i.e. the minimization of $\sigma[Area]$. The average and standard deviation across $50$ optimization repeats are shown.
    (d) Objective function values sampled during all optimization runs. The arrows indicate the typical trajectory of the optimizations, which first try to achieve values of $\mathbb{E}[Area]$ above 1000 a.u. and then try to minimize $\sigma[Area]$. A Pareto front that describes the trade-off between the two objectives becomes visible, as larger area's expectation values are accompanied by larger variability.
    (e) Objective function values sampled during a sample optimization run. Each experiment is color-coded (yellow to dark green) to indicate at which stage of the optimization it was performed. Exploration (white rim) and exploitation (black rim) points are indicated, as \textit{Gryffin} explicitly alternates between these two strategies. Later exploitation points (dark green, black rim) tend to focus on the minimization of $\sigma[Area]$, having already achieved $\mathbb{E}[Area] > 1000$ a.u.
    }
    \label{fig:hplc-opt}
\end{figure*}

As a realistic and challenging example, we consider the optimization of the aforementioned HPLC sampling protocol under the presence of significant noise in P1 and P3 (noisy query setting). In this first instance, we assume that the other conditions contain little noise and can thus be approximated as noiseless. As before, we consider normally distributed noise, truncated at zero. We assume a standard deviation of $0.008$ mL for P1, and $0.08$ mL for P3. In this example, we assume we are aware of the presence of input noise in these parameters, and are interested in achieving a protocol that returns an expected peak area, $\mathbb{E}[Area]$, of at least $1000$ a.u. As a secondary objective, we would like to minimize the output variability, $\sigma[Area]$, as much as possible while maintaining $\mathbb{E}[Area] > 1000$ a.u. 

To achieve the optimization goals, we use \golem to estimate both $\mathbb{E}[Area]$ and $\sigma[Area]$ as the optimization proceeds (Figure \ref{fig:hplc-opt}a). We then use \textit{Chimera}\cite{Hase:2018_chimera} to scalarize these two objectives into a single robust and multi-objective function, $g[Area]$, to be optimized. \textit{Chimera} is a scalarizing function that enables multi-objective optimization via the definition of a hierarchy of objectives and associated target values. As opposed to the post-hoc analysis discussed in the previous section, in this example we start with no prior experiment being available and let the optimization algorithm request new experiments in order to identify a suitable protocol. Here we perform virtual HPLC runs using \textit{Olympus}\cite{Hase:2020_olympus}, which allows to simulate experiments via Bayesian Neural Network models. These probabilistic models capture the stochastic nature of experiments, such that they return slightly different outcomes every time an experiment is simulated. In other words, they simulate the heteroscedastic noise present in the experimental measurements. While measurement noise is not the focus of this work, it is another source of uncertainty routinely encountered in an experimental setting. As such, it is included in this example application. Bayesian optimization algorithms are generally robust to some level of measurement noise, as this source of uncertainty is inferred by the surrogate model. However, the combination of output and input noise in the same experiment is particularly challenging, as both sources of noise manifest themselves as noisy measurements despite the different origin. In fact, in addition to measurement noise, here we inject input noise into the controllable parameters P1 and P3. Hence, while the optimization algorithm may request a specific value for P1 and P3, the actual, realized ones will differ. This setup therefore contains noise in both input experimental conditions and measurements.

While large input noise would be catastrophic in most standard optimization campaigns (as shown in section \ref{section:benchmarks_noisy}, Figure \ref{fig:xmin_examples}), \golem allows the optimization to proceed successfully. With the procedure depicted in Figure \ref{fig:hplc-opt}a, on average, \textit{Gryffin} was able to identify parameter settings that achieve $\mathbb{E}[Area] > 1000$ a.u. after less than $50$ experiments (Figure \ref{fig:hplc-opt}b). Equivalent results were obtained with \textit{GPyOpt} and \textit{Hyperopt} (Figure \ref{si_fig:hplc_opt_gpyopt_hopt}). The improvements in this objective are, however, accompanied by a degradation in the second objective, output variability, as measured by $\sigma[Area]$ (Figure \ref{fig:hplc-opt}c). This effect is due to the inevitable trade-off between the two competing objectives being optimized. After having reached its primary objective, the optimization algorithm mostly focused on improving the second objective, while satisfying the constraint defined for the first one. This behavior is visible in Figure \ref{fig:hplc-opt}d-e. Early in the optimization, \textit{Gryffin} is more likely to query parameter settings with low $\mathbb{E}[Area]$ and $\sigma[Area]$ values. At a later stage, with more information about the response surface, the algorithm focused on lowering $\sigma[Area]$ while keeping $\mathbb{E}[Area]$ above $1000$ a.u. Due to input uncertainty, the Pareto front highlights an irreducible amount of output variance for any non-zero values of expected area (Figure \ref{fig:hplc-opt}d). An analysis of the true robust objectives shows that, given the $\mathbb{E}[Area] > 1000$ a.u. constraint, the best achievable $\sigma[Area]$ values are $\sim 300$ a.u. (Figure \ref{si_fig:hplc_true_optimum}).

The traces showing the optimization progress (Figure \ref{fig:hplc-opt}b-c) display considerable spread around the average performance. This is expected and due to the fact that both $\mathbb{E}[Area]$ and $\sigma[Area]$ are estimates based on scarce data, as they cannot be directly observed. As a consequence, these estimates fluctuate as more data is collected. In addition, it may be the case that while \golem estimates $\mathbb{E}[Area]$ to be over $1000$ a.u., its true value for a certain set of input conditions may actually be below $1000$, and vice-versa. In fact, at the end of the $50$ repeated optimization runs, $10$ (i.e., $20$\%) of the identified optimal solutions had true $\mathbb{E}[Area]$ below $1000$ a.u. (this was the case for $24$\% of the optimizations with \textit{GPyOpt}, and $34$\% for those with \textit{Hyperopt}). However, when using ensemble trees as the surrogate model, it is possible to obtain an estimate of uncertainty for \golem's expectation estimates. With this uncertainty estimate, one can control the probability that \golem's estimates satisfy the objective's constraint that was set. For instance, to have a high probability of the estimate of $\mathbb{E}[Area]$ being above $1000$ a.u., we can setup the optimization objective in \textit{Chimera} with the constraint that $\mathbb{E}[Area] - 1.96 \times \sigma (\mathbb{E}[Area]) > 1000$ a.u., which corresponds to optimizing against the lower bound of the $95$\% confidence interval of \golem's estimate. Optimizations set up in this way correctly identified optimal solutions with $\mathbb{E}[Area] > 1000$ a.u. in all 50 repeated optimization runs (Figure \ref{si_fig:hplc_opt_lcb}).

As a final test, we simulate the example above, in which we targeted the optimization of the lower-bound estimate of $\mathbb{E}[Area]$, with all experimental conditions containing a considerable amount of noise. For all input variables we consider normally distributed noise truncated at zero, with a standard deviation of $0.008$ mL for P1, $0.06$ mL for P2, $0.08$ mL for P3, $0.2$ mL/min for P4, $8$ Hz for P5, and $1$ s for P6. This is an even more challenging optimization scenario, with input noise compounding from all variables. In this case, \textit{Hyperopt} achieved $\mathbb{E}[Area] > 1000$ a.u. after about $100$ experiments on average, \textit{Gryffin} achieved $\mathbb{E}[Area]$ values around the targeted value of $1000$ a.u. after $120-130$ experiments, and \textit{GPyOpt} only when close to $200$ experiments (Figure \ref{si_fig:hplc_opt_lcb_all_noisy}). As expected, the noisier the experimental conditions (larger noise and/or more noisy variables) the less efficient the optimization. However, \golem still enabled the algorithms tested to achieve the desired objective of $\mathbb{E}[Area] > 1000$ within the pre-defined experimental budget. After $200$ experiments, \textit{Hyperopt} correctly identified solutions with $\mathbb{E}[Area] > 1000$ a.u. in $78$\% of the optimization runs, \textit{GPyOpt} in $70$\%, and \textit{Gryffin} in $42$\%. We stress that \golem is not a substitute to developing precise experimental protocols. A noise-free (or reduced-noise) experimental protocol will always allow for faster optimization and better average performance. While \golem can mitigate the detrimental effects of input noise on optimization, it is still highly desirable to minimize noise in as many input conditions as possible.

This example application shows how \golem can easily be integrated into a Bayesian optimization loop for the optimization of experimental protocols with noisy experimental conditions.


\section{Conclusion}
In summary, \golem provides a simple, inexpensive, yet flexible approach for the optimization of experimental protocols under noisy experimental conditions. It can be applied retrospectively, for the analysis of previous results, as well as on-the-fly in conjunction with most experiment planning strategies to drive optimizations toward robust solutions. The approach was found to perform particularly well when used with systematic searches and Bayesian optimization algorithms. Optimization under noisy conditions is considerably more challenging than typical optimization tasks. When such noise is known but cannot be removed or corrected for, \golem enables optimizations that would otherwise be infeasible. An open-source implementation of the \golem algorithm is available on GitHub\cite{github_repo} under an MIT license.

	
	
	

	\section*{Acknowledgments}

		The authors thank Melodie Christensen and Lars Yunker for valuable and insightful discussions. The authors acknowledge the generous support of the National Research Council (NRC) of the Government of Canada. M.A. is supported by a Postdoctoral Fellowship of the Vector Institute. F.H. acknowledges financial support from the Herchel Smith Graduate Fellowship and the Jacques-Emile Dubois Student Dissertation Fellowship. R.J.H. gratefully acknowledges the Natural Sciences and Engineering Research Council of Canada (NSERC) for provision of the Postgraduate Scholarships-Doctoral Program (PGSD3- 534584-2019). A.A.-G. acknowledges support form the Canada 150 Research Chairs program and CIFAR. A.A.-G. acknowledges the generous support from Anders G. Fr\"{o}seth. All computations reported in this paper were completed on the computing clusters of the Vector Institute. Resources used in preparing this research were provided, in part, by the Province of Ontario, the Government of Canada through CIFAR, and companies sponsoring the Vector Institute. Finally, we thank the anonymous reviewers for their critical reading and comments that resulted in a much-improved manuscript.


	\phantomsection\addcontentsline{toc}{section}{\refname}\putbib[main]
\end{bibunit}

\clearpage
\newpage

\begin{bibunit}[unsrt]

	
	\setcounter{page}{1}
	\setcounter{section}{0}
	\setcounter{subsection}{0}
	\setcounter{figure}{0} 
	
        \renewcommand{\thesection}{S.\arabic{section}}
	\renewcommand{\thefigure}{S\arabic{figure}}
	\renewcommand{\thetable}{S\arabic{table}}
	
	
	\onecolumngrid
	\subsection*{\normalsize{Supplementary Information}{\\}{\vspace{6pt}}
	 			\large{Golem: An algorithm for robust experiment and process optimization}{\\}{\vspace{6pt}}
			        \normalsize{\normalfont{Matteo Aldeghi,$^{1,2,3}$ Florian H\"{a}se,$^{4,1,2,3}$ Riley J. Hickman,$^{2,3}$ Isaac Tamblyn,$^{1,5}$ Al\'{a}n Aspuru-Guzik$^{1,2,3,6}$}}
			        {\\}{\vspace{6pt}}
			        \small{\normalfont{$^1$\textit{Vector Institute for Artificial Intelligence, Toronto, ON, Canada}}}{\\}
			        \small{\normalfont{$^2$\textit{Chemical Physics Theory Group, Department of Chemistry, University of Toronto, Toronto, ON, Canada}}}{\\}
			        \small{\normalfont{$^3$\textit{Department of Computer Science, University of Toronto, Toronto, ON, Canada}}}{\\}
			        \small{\normalfont{$^4$\textit{Department of Chemistry and Chemical Biology, Harvard University, Cambridge, MA, USA}}}{\\}
			        \small{\normalfont{$^5$\textit{National Research Council of Canada, Ottawa, ON, Canada}}}{\\}
			        \small{\normalfont{$^6$\textit{Lebovic Fellow, Canadian Institute for Advanced Research, Toronto, ON, Canada}}}
			    }
	{\vspace{18pt}}

\section{Formulating Golem}
\label{section:golem_derivation}

Consider a sequential optimization in which the goal is to find a set of input conditions $\bm{x} \in \mathcal{X}$ corresponding to the global minimum of the function $f: \mathcal{X} \mapsto \mathbb{R}$,

\begin{align} 
\bm{x}^* = \argmin_{\bm{x} \in \mathcal{X}} f(\bm{x}).
\end{align}

At each iteration $k$, we query parameter values $\bm{x}_k \in \mathcal{X}$, where $\mathcal{X}$ is a compact subset of Euclidean space $\mathbb{R}^D$ and $D \in \mathbb{N}^*$. However, while the desired query location $\bm{x}$ can be precisely controlled, uncertainty in the execution results in $\widetilde{\bm{x}}=\bm{x} + \bm{\delta}$ being the input realized, where $\bm{\delta}$ is a random variable with probability density $p(\bm{\delta} | \bm{x})$. Thus, after $K$ optimization iterations, we will have built a dataset $\widetilde{\mathcal{D}}_K= \{ \bm{x}_k, f(\widetilde{\bm{x}}_k) \} _{k=1}^K$, where uncertainty in the inputs causes stochastic evaluations of the response function $f(\bm{x})$. A noiseless dataset $\mathcal{D}_K = \{ \bm{x}_k, f(\bm{x}_k) \} _{k=1}^K$ can be obtained when $p(\bm{\delta} | \bm{x})$ is a delta function.

We seek a robustness measure $g(\bm{x}_k)$ that will allow us to estimate the merit of each solution $\bm{x}_k$ given its uncertainty. This would not only allow for \textit{post-hoc} rescaling of the merit of each parameter location based on its robustness, but also to direct experiment planning algorithms towards robust optima during the optimization campaign.

A natural choice is to take the expectation of $f(\widetilde{\bm{x}}_k)$ by marginalizing over the input uncertainty, such that $g(\bm{x}_k) = \mathbb{E}[f(\widetilde{\bm{x}}_k)]$. Consider $\bm{x}_k$ to be the query location at iteration $k$, where we would like to evaluate $f(\bm{x}_k)$. The location actually realized, due to the uncertainty, is $\widetilde{\bm{x}}_k = \bm{x}_k + \bm{\delta}$, where $\bm{\delta}$ is a random variable with probability density $p(\bm{\delta} | \bm{x}_k)$. In principle, the expectation of $f(\widetilde{\bm{x}}_k)$ given $p(\widetilde{\bm{x}}_k)$ can be obtained as follows:

\begin{align} 
\mathbb{E}[f(\widetilde{\bm{x}}_k)] &= \int_{\mathbb{R}^D} f( \bm{x} ) p(\widetilde{\bm{x}}_k) d\bm{x}   \label{eq:integral} \\
 &= \int_{\mathbb{R}^D} f(\bm{x}) p(\bm{x}_k + \bm{\delta}) d\bm{x},
\end{align}

with integration over the support of $p(\widetilde{\bm{x}}_k)$, because $p(\bm{\delta} | \bm{x}_k)$ may extend beyond the bounds of $\mathcal{X}$. Here we assume integration over $\mathbb{R}^D$ without loss of generality. If physical bounds are present, this can be reflected in the support of $p(\bm{\delta}|\bm{x}_k)$; for instance, the uncertainty in volume of dispensed liquid may be modelled by a gamma distribution, given that volume dispensed can only be equal to or greater than zero. Eq. \ref{eq:integral} is another way to write the expectation obtained by marginalising over the noise factors $\bm{\delta}$,

\begin{align} 
\mathbb{E}[f(\widetilde{\bm{x}}_k)] = \int_{\mathbb{R}^D} f( \bm{x}_k + \bm{\delta}) p({\bm{\delta}}) d\bm{\delta}.
\end{align}

More explicitly,

\begin{align} 
\mathbb{E}[f(\widetilde{\bm{x}}_k)] &= \int_{\mathbb{R}^D} f( \bm{x}_k + \bm{t}) p_{\bm{\delta}}(\bm{t}) d\bm{t} .
\end{align}

Defining $\bm{x} = \bm{x}_k + \bm{t}$, $d\bm{x} = d\bm{t}$, and given that $p_{\bm{\delta}}(\bm{x} - \bm{x}_k) = p_{\bm{x}_k + \bm{\delta}}(\bm{x})$,

\begin{align} 
\mathbb{E}[f(\widetilde{\bm{x}}_k)] &= \int_{\mathbb{R}^D} f(\bm{x}) p_{\bm{\delta}}(\bm{x} - \bm{x}_k) d\bm{x} \\
&= \int_{\mathbb{R}^D} f(\bm{x}) p_{\bm{x}_k + \bm{\delta}}(\bm{x}) d\bm{x}.
\end{align}

Eq. \ref{eq:integral} as shown above is recovered when simplifying the notation with $p_{\bm{x}_k + \bm{\delta}}(\bm{x}) = p(\bm{x}_k + \bm{\delta})$. Finally, note that other integrated measures of robustness can be defined, for instance considering higher moments of the objective function.\cite{Beyer:2007,Beland:2017}

While $f(\bm{x})$ is unknown, an approximation $\hat{f}(\bm{x})$ can be built from $\mathcal{D}_K$ or $\widetilde{\mathcal{D}}_K$. For simplicity, we assume $\hat{f}( \bm{x} ) \approx f(\bm{x})$ and from now on will refer to $\hat{f}( \bm{x} )$ simply as $f(\bm{x})$. If we can solve the above integral efficiently, we can then use $\mathbb{E}[f(\widetilde{\bm{x}}_k)]$ as the robust merit $g(\bm{x})$ for each condition $\bm{x}_k$, and $\mathcal{G}_K= \{ \bm{x}_k, g_k \} _{k=1}^K$ could be used with an experiment planning algorithm of choice to solve the robust optimization problem as

\begin{align} 
\bm{x}^* = \argmin_{\bm{x} \in \mathcal{X}} g(\bm{x}).
\end{align}

However, there is no closed form solution of Eq. \ref{eq:integral} for most combinations of $f(\bm{x})$ and $p(\widetilde{\bm{x}}_k)$. Its numerical approximation is expensive, becoming intractable with increasing dimensionality and number of samples.

\subsection{Continuous input variables}

The space over which $p(\widetilde{\bm{x}})$ is supported (here considered to be $\mathbb{R}^D \supset \mathcal{X}$) can be partitioned into $M \in \mathbb{N}^*$ non-overlapping tiles $\{\mathcal{T}^D_m\}_{m=1}^{M}$, with $\mathcal{T}^D_m \subset \mathbb{R}^D\ \forall\ m \leq M$, to create a $D$-dimensional tessellation. With this discretization, Eq. \ref{eq:integral} can be decomposed into a finite series with integration over each tile $\mathcal{T}^D_m$:

\begin{align} 
\mathbb{E}[f(\widetilde{\bm{x}}_k)] = \sum_{m=1}^{M} \int_{\mathcal{T}_m^D} f(\bm{x}) p(\widetilde{\bm{x}}_k) d\bm{x}.
\end{align}

Assuming a piecewise constant model of $f(\bm{x})$, such as a regression tree, $f(\bm{x})$ is constant within the partition $\mathcal{T}_m^D$ and can be brought outside the integral:

\begin{align} 
\mathbb{E}[f(\widetilde{\bm{x}}_k)] = \sum_{m=1}^{M} f_m \int_{\mathcal{T}_m^D} p(\widetilde{\bm{x}}_k) d\bm{x}.
\label{eq:tiles}
\end{align}

The integral over $\mathcal{T}_m^D$

\begin{align} 
\int_{\mathcal{T}_m^D} p(\widetilde{\bm{x}}_k) d\bm{x} &= P(\widetilde{\bm{x}}_k \in \mathcal{T}_m^D)
\end{align}

is the probability of $\bm{x}_k$ being in tile $m$, given the uncertainty $p(\bm{\delta}_k)$. Therefore, $\mathbb{E}[f(\widetilde{\bm{x}}_k)]$ is effectively a weighted average,  

\begin{align} 
\mathbb{E}[f(\widetilde{\bm{x}}_k)] = \sum_{m=1}^{M} f_m \cdot P(\widetilde{\bm{x}}_k \in \mathcal{T}_m^D),
\end{align}

where all possible outcomes are weighted by their probability given the targeted parameter location $\bm{x}_k$. Assuming independent input uncertainties, $P(\widetilde{\bm{x}}_k \in \mathcal{T}_m^D)$ can be factorized:

\begin{align} 
P(\widetilde{\bm{x}}_k \in \mathcal{T}_m^D) &= \prod_{d=1}^D P(\widetilde{x}_{k,d} \in \mathcal{T}_{m, d}).
\end{align}

The probability $P(\widetilde{x}_{k,d} \in \mathcal{T}_{m, d})$ is obtained from the cumulative distribution function $F_{k,d}$ of $p(\widetilde{\bm{x}}_k)$, evaluated at the upper and lower bounds of tile $m$ in dimension $d$,

\begin{align} 
P(\widetilde{x}_{k,d} \in \mathcal{T}_{m,d}) &=  F_{k,d}(\max_d\ \mathcal{T}_m^D) - F_{k,d}(\min_d\ \mathcal{T}_m^D),
\label{eq:cdf}
\end{align}

where $\max_d\ \mathcal{T}_m^D$ and $\min_d\ \mathcal{T}_m^D$ are the upper and lower bounds of tile $m$, respectively, in the $d$ dimension. $F_{k,d}(a) = \int_{-\infty}^{a} p(\widetilde{x}_{k,d}) dx_d$ is the cumulative distribution function of $p(\widetilde{x}_k)$ in the $d$ dimension. It thus follows that, combining Eq. \ref{eq:tiles}-\ref{eq:cdf}, for a piecewise constant model $f(\bm{x})$ and any parametric distribution $p(\widetilde{\bm{x}}_k)$ with known $F_k$, the desired expectation can be computed as

\begin{align} 
\mathbb{E}[f(\widetilde{\bm{x}}_k)] = \sum_{m=1}^{M} f_m \prod_{d=1}^D [ F_{k,d}(\max_d\ \mathcal{T}_m^D) - F_{k,d}(\min_d\ \mathcal{T}_m^D) ].
\label{eq:golem}
\end{align}

In this work, we model $f(\bm{x})$ with single regression trees as well as their ensemble variants, like random forest and extremely randomized trees.\cite{Breiman:2001,Geurts:2006} When ensembles are used, a different expectation $\mathbb{E}[f_t(\widetilde{\bm{x}}_k)]$ is obtained for each tree $t$, in which case we take their average as the most reliable estimate of robustness:

\begin{align} 
\mathbb{E}[f(\widetilde{\bm{x}}_k)] = \frac{1}{T} \sum_{t=1}^{T} \sum_{m=1}^{M} f_{t,m} \prod_{d=1}^D [ F_{k,d}(\max_d\ \mathcal{T}_{t,m}^D) - F_{k,d}(\min_d\ \mathcal{T}_{t,m}^D) ].
\label{eq:full_golem}
\end{align}

\subsection{Discrete input variables}
Tree-based machine learning approaches can also take discrete and categorical variables as features, such that uncertainty in these types of inputs can also be handled by \golem. When optimizing over a discrete space $\mathcal{X} \subset \mathbb{N}^D$, both $\bm{x} \in \mathcal{X}$ and $\bm{\delta} \in \mathbb{N}^D$ are discrete and the expectation of $f$ is expressed as a sum over $\mathbb{N}^D$,

\begin{align} 
\mathbb{E}[f(\widetilde{\bm{x}}_k)] &= \sum_{\mathbb{N}^D} f( \bm{x} ) p(\widetilde{\bm{x}}_k),
\end{align}

where $p(\widetilde{\bm{x}}_k)$ is a discrete probability distribution. The $\mathbb{N}^D$ space supporting this distribution can be partitioned into $M \in \mathbb{N}^*$ non-overlapping tiles $\{\mathcal{T}^{D}_m\}_{m=1}^{M}$, with $\mathcal{T}^{D}_m \subset \mathbb{N}^{D}\ \forall\ m \leq M$,

\begin{align} 
\mathbb{E}[f(\widetilde{\bm{x}}_k)] &= \sum_{m=1}^{M}  \sum_{\mathcal{T}_m^{D}} f(\bm{x}) p(\widetilde{\bm{x}}_k).
\end{align}

For a model that takes a constant $f_m$ value within each tile $\mathcal{T}_m^{D}$,

\begin{align} 
\mathbb{E}[f(\widetilde{\bm{x}}_k)] &=  \sum_{m=1}^{M} \Big( f_m \sum_{\mathcal{T}_m^{D}} p(\widetilde{\bm{x}}_k) \Big) \\
&= \sum_{m=1}^{M} f_m \cdot P(\widetilde{\bm{x}}_k \in \mathcal{T}_m^{D}).
\end{align}

As for continuous variables, we assume independent uncertainty across input variables, such that $P(\widetilde{\bm{x}}_k \in \mathcal{T}_m^{D})$ factorises. The rest of the derivation then follows the same argument as for continuous variables. Thus, the only difference to continuous variables is that $\widetilde{\bm{x}}$ and its probability distributions are discrete.

\subsection{Categorical input variables}
If the optimization occurs over a $D$-dimensional space with $C \in \mathbb{N}^*$ categorical options, at each iteration we query a point $\bm{x}_k \in \mathbb{S}^{D\times C}$ that selects a category $\bm{z}_d$ for each dimension $d$. This information can be encoded as $C$-dimensional one-hot encoded vectors, $\mathbb{S}^{D\times C} = \{\bm{z} \in \mathbb{R}^{D\times C} | z_{d,c} \in \{0,1\}; \sum_{c=1}^C z_{d,c} = 1\ \forall\ d \leq D \}$. The uncertainty over categorical variables can then be represented by any suitable probability distribution on the simplex, $p(\widetilde{\bm{x}}_{k}) \in \Delta^{D\times (C-1)} = \{ p(\bm{z}) \in \mathbb{R}^{D\times C} | p(z_{d,c}) \in [0,1]; \sum_{c=1}^C p(z_{d,c}) = 1\ \forall\ d \leq D\}$. In this scenario, the expectation of $f$ queried at location $\bm{x}_k$ and considering the uncertainty due to $\bm{\delta}_k$ is

\begin{align} 
\mathbb{E}[f(\widetilde{\bm{x}}_k)] &= \sum_{\mathbb{S}^{D\times C}} f( \bm{x} ) p(\widetilde{\bm{x}}_k).
\end{align}

Similar to what was done before, we partition the space $\mathbb{S}^{D\times C}$ in $M \in \mathbb{N}^*$ non-overlapping tiles $\{\mathcal{T}^{D\times C}_m\}_{m=1}^{M}$, with $\mathcal{T}^{D\times C}_m \subset \mathbb{S}^{D\times C}\ \forall\ m \leq M$:

\begin{align} 
\mathbb{E}[f(\widetilde{\bm{x}}_k)] &=  \sum_{m=1}^{M}  \sum_{\mathcal{T}_m^{D\times C}} f( \bm{x} ) p(\widetilde{\bm{x}}_k).
\end{align}

For a model that assumes a constant $f_m$ values within each tile $\mathcal{T}_m^{D\times C}$,

\begin{align} 
\mathbb{E}[f(\widetilde{\bm{x}}_k)] &=  \sum_{m=1}^{M} \Big( f_m \sum_{\mathcal{T}_m^{D\times C}} p(\widetilde{\bm{x}}_k) \Big) \\
&= \sum_{m=1}^{M} f_m \cdot P(\widetilde{\bm{x}}_k \in \mathcal{T}_m^{D\times C}).
\end{align}

Assuming independent uncertainty across input variables,

\begin{align} 
P(\widetilde{\bm{x}}_k \in \mathcal{T}_m^{D\times C}) = \prod_{d=1}^D P(\widetilde{\bm{x}}_{k,d} \in \mathcal{T}_{m, d}^{1\times C}).
\end{align}

The probability of $\widetilde{\bm{x}}_{k,d}$ being in the $1\times C$ dimensional tile $\mathcal{T}_{m, d}^{1\times C}$ can be readily computed from the user-defined probabilities $p(z_{d,c})$ indicating the uncertainty over categorical variables, such that

\begin{align} 
P(\widetilde{\bm{x}}_{k,d} \in \mathcal{T}_{m, d}^{1\times C}) &= \sum_{c=1}^C  p(z_{d,c}) \cdot \mathbb{I}(\bm{z}_d \in \mathcal{T}_{m,d}^{1\times C}),
\end{align}

where $\mathbb{I}$ is the indicator function, taking the value of $1$ if the category $\bm{z}_d$ is in tile $\mathcal{T}_{m,d}^{1\times C}$ and $0$ otherwise. In our implementation, we build trees until leaves are pure, which means that each tile $\mathcal{T}_{m, d}^{1\times C}$ will contain a single category (when at least one sample per category is present). However, this does not necessarily need to be the case and one might decide to limit tree depth for computational efficiency.

\subsection{Multi-objective optimization}
\label{section:si_multiobj}

\begin{figure*}[htb]
    \centering
    \includegraphics[width=1.0\textwidth]{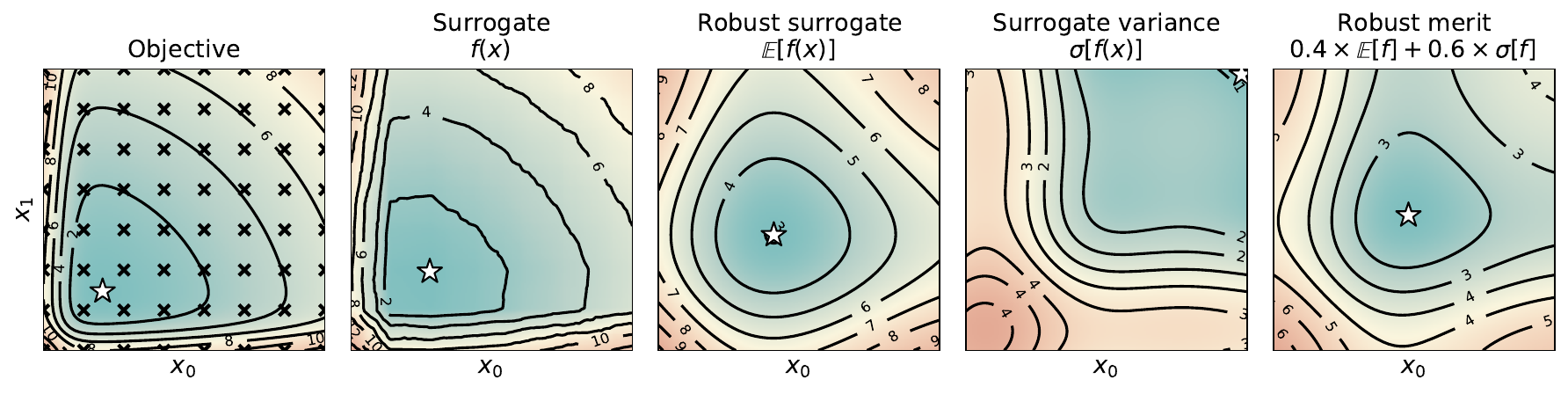}
    \caption{Multi-objective optimization with \golem. The first plot on the left shows the \textit{Cliff} objective function (section \ref{section:si_surfaces}). \golem's surrogate model was built using the $64$ samples marked as black crosses. The model, in this case, was a forest of $100$ extremely randomized trees\cite{Geurts:2006}. The robust surrogate was built assuming normally distributed input noise, in both dimensions, and with unit standard deviation. The variability of the objective function under this noise model was computed as $\sigma [f(\widetilde{\bm{x}})] = (\mathbb{E}[f(\widetilde{\bm{x}})^2] - \mathbb{E}[f(\widetilde{\bm{x}})]^2))^{\frac{1}{2}}$. The two objectives were combined into a single function to be optimized via the weighted sum $g(\bm{x}) = 0.4 \times \mathbb{E}[f(\widetilde{\bm{x}})] + 0.6 \times \sigma [f(\widetilde{\bm{x}})]$, where $0.4$ and $0.6$ are user-defined coefficients.}
    \label{si_fig:multiobj}
\end{figure*}

In addition to computing the expectation of $f(\bm{x})$, one can also consider its variance. As lower variance favors reproducibility, one might be interested in minimizing both $\mathbb{E}[f(\widetilde{\bm{x}}_k)]$ and $\sigma [f(\widetilde{\bm{x}}_k)] = Var[f(\widetilde{\bm{x}}_k)]^{\frac{1}{2}}$. The variance can easily be obtained as $Var[f(\widetilde{\bm{x}}_k)] = \mathbb{E}[f(\widetilde{\bm{x}}_k)^2] - \mathbb{E}[f(\widetilde{\bm{x}}_k)]^2$. With both $\mathbb{E}[f(\widetilde{\bm{x}}_k)]$ and $\sigma [f(\widetilde{\bm{x}}_k)]$ available, one can carry out a multi-objective optimization by building a robust merit $g(\bm{x})$ that takes both objectives into account, using any scalarizing function of choice. Figure \ref{si_fig:multiobj} shows an example where a robust merit function is built via a weighted sum.

\subsection{\golem's assumptions}
\golem relies on three fundamental assumptions for its derivation as well as successful deployment. First, it is assumed that a piece-wise constant model, such as those obtained with tree-based algorithms, is able to provide an accurate surrogate model of the underlying objective function $f(\bm{x})$. In addition, it is assumed that such a surrogate model can be built from finite datasets $\mathcal{D}_K$ or $\widetilde{\mathcal{D}}_K$. It is expected that building an accurate surrogate model will be more challenging when using a noisy dataset $\widetilde{\mathcal{D}}_K$ based on stochastic queries of input conditions, than when using a noiseless dataset $\mathcal{D}_K$ based on deterministic ones. Second, it is assumed that the user knows and is able to accurately model input uncertainty via a parametric probability distribution $p({\bm{x}})$. Finally, a necessary assumption in the above \golem's derivation is that the uncertainties of different input conditions are independent, such that, for instance, $p(x_i,x_j) = p(x_i)p(x_j)$ where $i$ and $j$ are two different input dimensions. This assumption might not always be satisfied depending on the input conditions and experimental setup. For instance, imagine that some uncertainty is associated with both target temperature and dispensed volume of liquid for a hypothetical experiment. If the liquid is first dispensed at room temperature and then heated to the desired target temperature, the errors in volume dispensed and target temperature are indeed likely independent. However, if the liquid is heated to a target temperature before dispensing, the (unknown, realized) temperature might affect viscosity, which in turn will have an effect on dispensing errors. That said, \golem allows for the uncertainty in one input variable to depend on the query location of all input variables, such that one can define $p(\widetilde{x}_i \mid x_i, x_j)$. Using the same example from above, it is thus possible to specify how the dispensing uncertainty depends upon the target temperature, though not the realized, unknown one.

\subsection{Computational scaling}
\begin{figure*}[htb]
    \centering
    \includegraphics[width=1.0\textwidth]{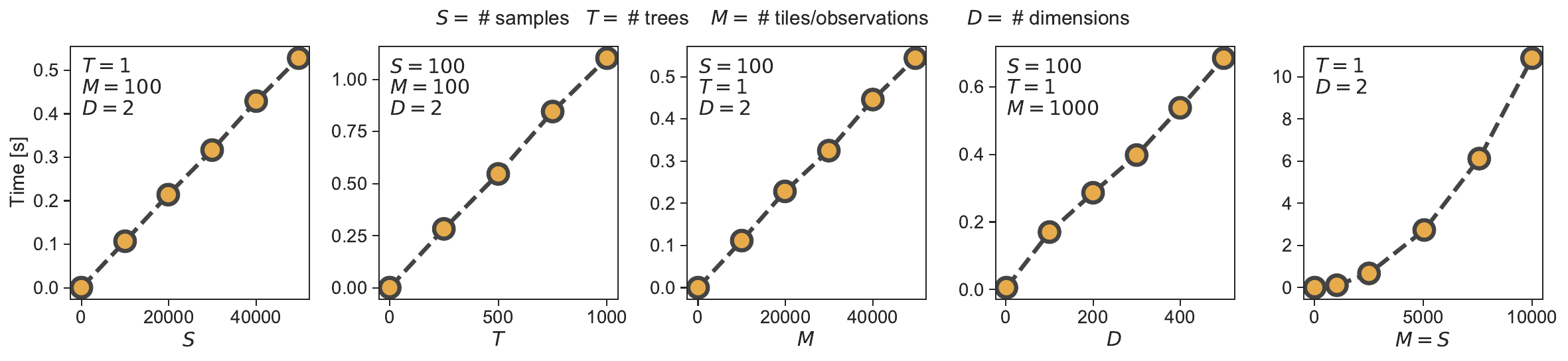}
    \caption{Computational scaling of \golem with respect to the number of predicted samples (S), number of trees used as part of the surrogate model (T), number of leaves in each tree, which, in our implementation correspond to the number observations in the dataset used for training (M), and dimensionality of the optimization task (D).}
    \label{si_fig:scaling}
\end{figure*}

To obtain the estimate of the robust merit for an input location $\bm{x}$, \golem evaluates Eq. \ref{eq:full_golem} after having fitted the tree-based surrogate model. For $S$ input locations, this involves performing operations over all input dimensions $D$, number of tiles $M$, and number of trees $T$. The time complexity of the algorithm (given an already trained tree-based model) thus scales linearly with respect to all these variables, $\mathcal{O}(S \times T \times M \times D)$. If the trees are allowed to grow until each leaf contains a single observation, as done in our \golem implementation, the number of tiles $M$ corresponds to the number of observations $K$ in the dataset $\mathcal{D}_K$. Typically, these are the observations for which one would like to re-evaluate the merits. In addition, we expect the number of trees $T$ and the input dimensionality $D$ to generally be small with respect to the number of observations. Hence, in a typical asymptotic scenario we have that $M = S \gg T,D$ with \golem displaying a quadratic runtime $\mathcal{O}(n^2)$ that depends on the number of observations collected. The time complexity can, however, be further reduced to $\mathcal{O}(S)$ by defining a maximum tree depth that would bound $M$. Figure \ref{si_fig:scaling} shows \golem's run time when varying $S$, $T$, $M$, and $D$ as discussed above. Note that, despite the quadratic scaling of the implementation, the run time has a small prefactor. In practice, robust merits for thousands of samples can be estimated on a single CPU core in a matter of seconds. For instance, evaluating the robust merits of 2500 samples using a surrogate model with 10 trees, for a two-dimensional problem, takes approximately 7 seconds on a single core of a 1.4 GHz Quad-Core Intel Core i5-8257U processor.

\section{Synthetic benchmarks}

In the following, information regarding the synthetic benchmark functions used to evaluate \golem's performance, additional analyses, and the details of all results obtained, are provided.

\subsection{Benchmark functions}
\label{section:si_surfaces}

In this work, we use three objective functions, which, given different assumed input uncertainties, create the different robust objective functions used as synthetic benchmarks. We refer to these three objective functions as \textit{Bertsimas}, \textit{Cliff}, and \textit{Sine}. The \textit{Bertsimas} function is taken from previous work on robust optimization by Bertsimas \textit{et al.}\cite{Bertsimas2009}. It corresponds to the following nonconvex polynomial function for $x \in [-1,3.2]$ and $y \in [-0.5, 4.4]$:

\begin{align} 
f(x,y) = &\ 2x^6 - 12.2x^5 + 21.2x^4 + 6.2x - 6.4x^3 - 4.7x^2 +y^6 - 11y^5 + 43.3y^4 - 10y - 74.8y^3 \nonumber \\
& + 56.9y^2 - 4.1xy - 0.1y^2x^2 + 0.4y^2x + 0.4x^2y.
\end{align}

In this work we place an upper bound to the function codomain, such that the \textit{Bertsimas} function used in practice is $\min(f(x,y), 80)$. This was done to avoid the extremely large values present outside its optimization domain. The \textit{Cliff} function is introduced in this work and is defined as follows, with $\bm{x} \in [0, 5]^D$, where $D$ is the number of dimensions:

\begin{align} 
f(\bm{x}) = \sum_{d=1}^D \frac{10}{1 + 0.3e^{6x_d}} + 0.2x_d^2 .
\end{align}

The \textit{Sine} function is also introduced in this work and is defined as follows, with $\bm{x} \in [-1, 1]^D$:

\begin{align} 
f(\bm{x}) = \sum_{d=1}^D \sin(2 \pi x_d^2) + x_d^2 + 0.2x_d .
\end{align}

The global minimum of the \textit{Bertsimas} function is at $(x^*, y^*) = (2.8, 4.0)$, the minimum of \textit{Cliff} is at $\bm{x}^* = (1.02874)^D$, and that of \textit{Sine} is at $\bm{x}^* = (-0.85297)^D$.

Discrete versions of \textit{Bertsimas} and \textit{Cliff} were obtained by discretizing and scaling their domain onto a $22 \times 22$ grid, such that $x \in \mathbb{N} \mid 1 \leq x \leq 22$ and $y \in \mathbb{N} \mid 1 \leq y \leq 22$. In these cases, the gobal minima are found at $(x^*, y^*) = (20, 20)$ and $(x^*, y^*) = (5, 5)$ for the \textit{Discrete Bertsimas} and \textit{Discrete Cliff} functions, respectively.

The robust objective functions S1--S8 are obtained by transforming the above functions based on specific input distributions, as shown in Figure~\ref{fig:grids} and detailed in Table~\ref{tab:surfaces}. For continuous functions, we used the \textit{Normal}, \textit{Gamma}, and \textit{Uniform} distributions with various scales. For discrete functions, we used the \textit{Poisson} and \textit{Discrete Laplace}\cite{Inusah:2006} distributions. However, note that any parametric distribution can in principle be used to model input uncertainty. In our \golem package, we implemented a Gamma distribution parametrized by its standard deviation and with variable lower or upper bounds. Similarly, we allow shifting the Poisson distribution such that any lower bound can be chosen. Details of these implementations can be found in \golem's GitHub repository\cite{github_repo}.

For all objective functions in Table \ref{tab:surfaces}, a close numerical approximation of their corresponding robust objective was obtained with \golem using a dense grid of $40,000$ samples. These samples extended beyond the optimization domain of each objective function, to be able to accurately model the objective function across all accessible regions of input space. For surfaces associated with unbounded probability distributions, samples were taken up to two standard deviations away from the optimization domain boundaries.

\begin{table*}[htb]
    \begin{center}  
        \caption{Details of the synthetic benchmark functions used to evaluate \golem. $^\dagger$One standard deviation for Normal, Gamma, and Discrete Laplace distributions; range for Uniform distributions; not applicable to Poisson distributions as not parametrized by scale. $^\ddagger$ Measure of the improvement in robustness between the minimum of the objective function and that of the robust objective function, relative to the range of the co-domain of the robust objective function.}
        \label{tab:surfaces}
        \begin{tabular}{l l l l l l l}
        \toprule
            Label & Function & Optimization domain & Probability distribution & Scale$^\dagger$ & Support & Improvement$^\ddagger$ \\
        \hline
            S1 & Cliff & $x_i \in [0,5]$ & Normal & 1.0 & $x_i \in \mathbb{R}$ & $25$\%  \\
            S2 & Cliff & $x_i \in [0,5]$ & Gamma & 2.0 & $x_i \in (-\infty, 5]$ & $51$\% \\
            S3 & Bertsimas & \makecell[tl]{$x_0 \in [-1,3.2]$ \\ $ x_1 \in [-0.5,4.4]$} & Uniform & 1.5 & \makecell[tl]{$x_0 \in [-1.75,3.95]$ \\ $ x_1 \in [-1.25,5.15]$} & $34$\% \\
            S4 & Bertsimas & \makecell[tl]{$x_0 \in [-1,3.2]$ \\ $x_1 \in [-0.5,4.4]$} & Normal & 0.8 & $x_i \in \mathbb{R}$ & $53$\% \\
            S5 & Sine & $x_i \in [-1,1]$ & Uniform &  0.5 & $x_i \in [-1.25,1.25]$ & $26$\% \\
            S6 & Sine & $x_i \in [-1,1]$ & Normal & 0.2 & $x_i \in \mathbb{R}$ & $26$\% \\
            S7 & Discrete Cliff & $x_i \in \mathbb{N} \mid 1 \leq x \leq 22 $ & Discrete Laplace & 3 & $x_i \in \mathbb{N}$ & $18$\% \\
            S8 & Discrete Bertsimas & $x_i \in \mathbb{N} \mid 1 \leq x \leq 22 $ & Poisson & n.a. & $x_i \in \mathbb{N} \mid x \geq 1 $ & $74$\% \\
        \botrule
        \end{tabular}
    \end{center} 
\end{table*}

\subsection{Bias due to boundary effects}
\label{section:si_infsampling}
\begin{figure*}[htb]
    \centering
    \includegraphics[width=1.0\textwidth]{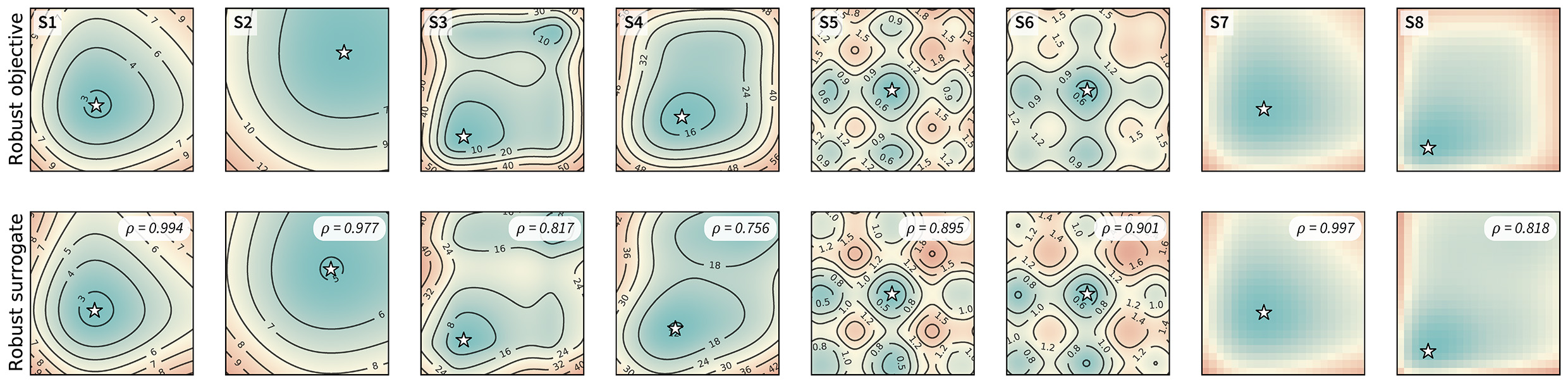}
    \caption{Converged robust surrogate estimates. The first row shows the true robust objective functions, which depend on the behavior of the objective function also outside the optimization domain shown. The second row shows converged robust surrogate estimates based on a regular grid of $10,000$ points within the optimization domain. Spearman's correlation ($\rho$) between the robust objective and its surrogate model are shown for each benchmark surface. Deviations from the ideal $\rho=1$ correlation are due to the inability of \golem's surrogate model (in this case, a single regression tree) to accurately capture the objective function's behavior beyond the optimization domain due to a lack of data in those regions.}
    \label{si_fig:infsampling}
\end{figure*}

\begin{figure*}[htb]
    \centering
    \includegraphics[width=1.0\textwidth]{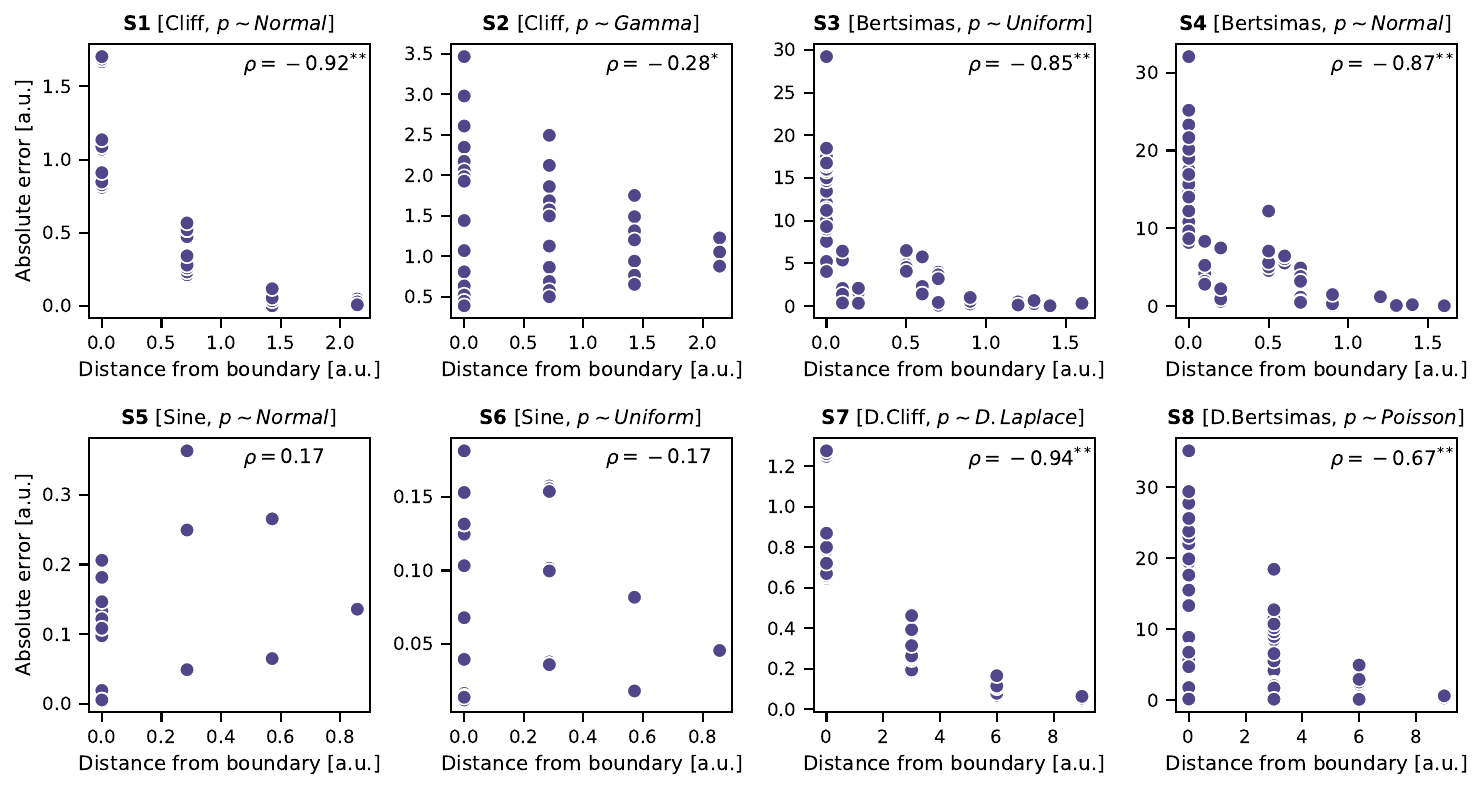}
    \caption{Relationship between distance from optimization boundaries and robustness estimate error. The data shown are for $64$ points uniformly sampled on a grid, as shown in the fourth row of Figure \ref{fig:grids}. The errors are the difference between the robustness estimates obtained with \golem based on $64$ datapoints and the ground truth estimate obtained as described in section \ref{section:si_surfaces}. On each plot we report the Spearman's correlation ($\rho$) between the errors and distances. Correlations with p-values less than $0.05$ are marked with one star, and those less than $0.01$ are marked with two. In the majority (six out of eight) of the test surfaces considered, there is a significant negative correlation between boundary distance and absolute errors.}
    \label{si_fig:boundary_error}
\end{figure*}

When the input parameters are noisy, the realized location of the queries does not correspond to that of the requested location. As a consequence, while one requests only locations within the bounds of the defined optimization domain, objective function evaluations outside of these bounds are possible. To know the true robustness of each solution within the optimization domain, one would thus need to know how the objective function behaves outside of the bounds of the optimization. As we approximate the objective function with a machine learning model, based on a dataset that has no samples outside the optimization domain, the surrogate model built is likely to be poor far outside the boundaries of the optimization domain. This lack of information results in a biased robust surrogate model also in the limit of infinite sampling within the optimization domain. This effect is exemplified by Figure \ref{si_fig:infsampling}, in which a surrogate robust objective was built with \golem using a dataset containing $10,000$ samples equally spaced within the optimization domain only. Another consequence of this boundary effect is that \golem's estimates of the robust objective tend to be less accurate for points close to the boundaries of the optimization domain (Figure~\ref{si_fig:boundary_error}).

\subsection{Cumulative robust regret as performance measure}
\label{section:si_regret}

To compare the relative optimization performance of all algorithms tested on specific benchmark functions $f(\cdot)$ we used the following definition of cumulative robust regret:

\begin{align}
\sum_{k=1}^{K} g\Big( \argmin_{\bm{x} \in \bm{x}_{1:k}} \hat{g}(\bm{x}_{1:k})\Big),
\end{align}

where $\bm{x}_{1:k}$ are all samples collected up to iteration $k$, $\hat{g}(\cdot)$ is the estimate of the robust merits obtained with \golem after training on $\mathcal{D}_{1:k}$, and $g(\cdot)$ is the true robust objective function. The true robust objective is obtained as an accurate \golem estimate by using a dense grid of $40,000$ samples, as mentioned in section \ref{section:si_surfaces}. To measure the performance of each optimization algorithm without \golem as the baseline, the original values of the merits, as obtained from the objective function $f(\cdot)$, were used instead of those derived with $\hat{g}(\cdot)$. Effectively, at each iteration, after having collected one additional sample, we estimate which one (among all samples collected thus far) has the best robust merit as estimated by \golem. Then, we take the true robust merit of this sample. All true robust merits obtained in this way for $k=1$ to $k=K$, where $K$ is the total number of samples, are then summed. This measure quantifies the speed at which the optimization algorithm has discovered better robust solutions. Given we are performing minimizations, the lower the cumulative regret, the better performing the algorithm is. For visualization and interpretation purposes, we normalize the values of cumulative regrets obtained in this way for each test on a specific benchmark surface. Hence, within each plot (e.g., in Figure~\ref{fig:noiseless_benchmarks}), a value of zero corresponds to the best cumulative regret observed for that surface, and a value of one to the worst. Note that, because this measure is not normalized across benchmark surfaces, comparisons are meaningful only with respect to a specific surface.

\subsection{Impact of uniformity and sampling of the boundaries}
\label{section:si_lowdisc}

\begin{figure*}[htb]
    \centering
    \includegraphics[width=1.0\textwidth]{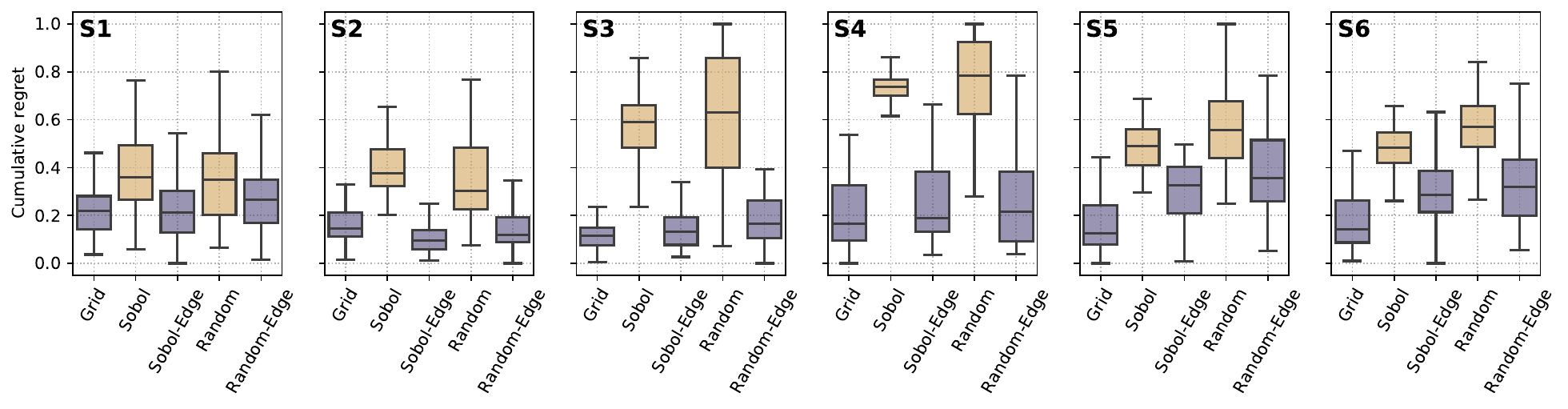}
    \caption{Optimization performance of different design of experiment approaches when used with \golem. Box plots show the distributions of cumulative regrets obtained across 50 optimization repeats. The boxes show the first, second, and third quartiles of the data, with whiskers extending up to 1.5 times the interquartile range. Boxes for approaches that sampled the boundaries of the optimization domain are shown in purple, while those for approaches that did not sample boundaries exactly are shown in yellow.}
    \label{si_fig:low-disc-seq}
\end{figure*}

In the optimization benchmarks carried out we noted that \textit{Grid} generally performed better than \textit{Random}. We hypothesized this might be caused by one or two features of these approaches. First, the difference in performance could be due the more uniform sampling of input space that is guaranteed with \textit{Grid}. To test this hypothesis, we performed optimizations, in the noiseless setting, with a Sobol sequence (we refer to this approach as \textit{Sobol}), which samples input space more uniformly than \textit{Random} but less than \textit{Grid}. Second, the performance difference could be due to the fact \textit{Grid} guarantees good sampling of the boundaries of the optimization domain. As discussed in section \ref{section:si_infsampling}, boundaries effect are present due to the lack of information on the objective function's behavior outside the optimization domain. To test this hypothesis, we benchmarked two additional approaches, in which we augmented \textit{Random} and \textit{Sobol} with samples at exactly the optimization domain boundaries. We placed these samples at the same locations of those in \textit{Grid}. We refer to these two approaches as \textit{Sobol-Edge} and \textit{Random-Edge}. In all cases, we allowed 196 objective function evaluations in total, as in the benchmarks described in section \ref{section:benchmarks}. Figure \ref{si_fig:low-disc-seq} summarizes the results obtained with the above approaches against benchmark functions S1--S6. We found that the unlikely sampling of the boundaries was the primary factor negatively impacting the performance of \textit{Random} when used in conjunction with \golem. The sampling of the boundaries present in \textit{Sobol-Edge} and \textit{Random-Edge} allowed \golem to better estimate the robustness of the input parameters collected, resulting in better performance for these two approaches as compared to \textit{Sobol} and \textit{Random}. The effect of uniform sampling was not noticeable in S1--S4, with \textit{Grid}, \textit{Sobol-Edge}, and \textit{Random-Edge} performing equally well. A small difference in performance was noticeable only for the rougher surfaces S5 and S6. There, \textit{Grid} performed better than \textit{Sobol-Edge}, which in turn was better than \textit{Random-Edge}, as expected if sampling uniformity were beneficial to robust optimizations with \golem.

\subsection{Influence of approximate surrogate model}
\begin{figure*}[htb]
    \centering
    \includegraphics[width=1.0\textwidth]{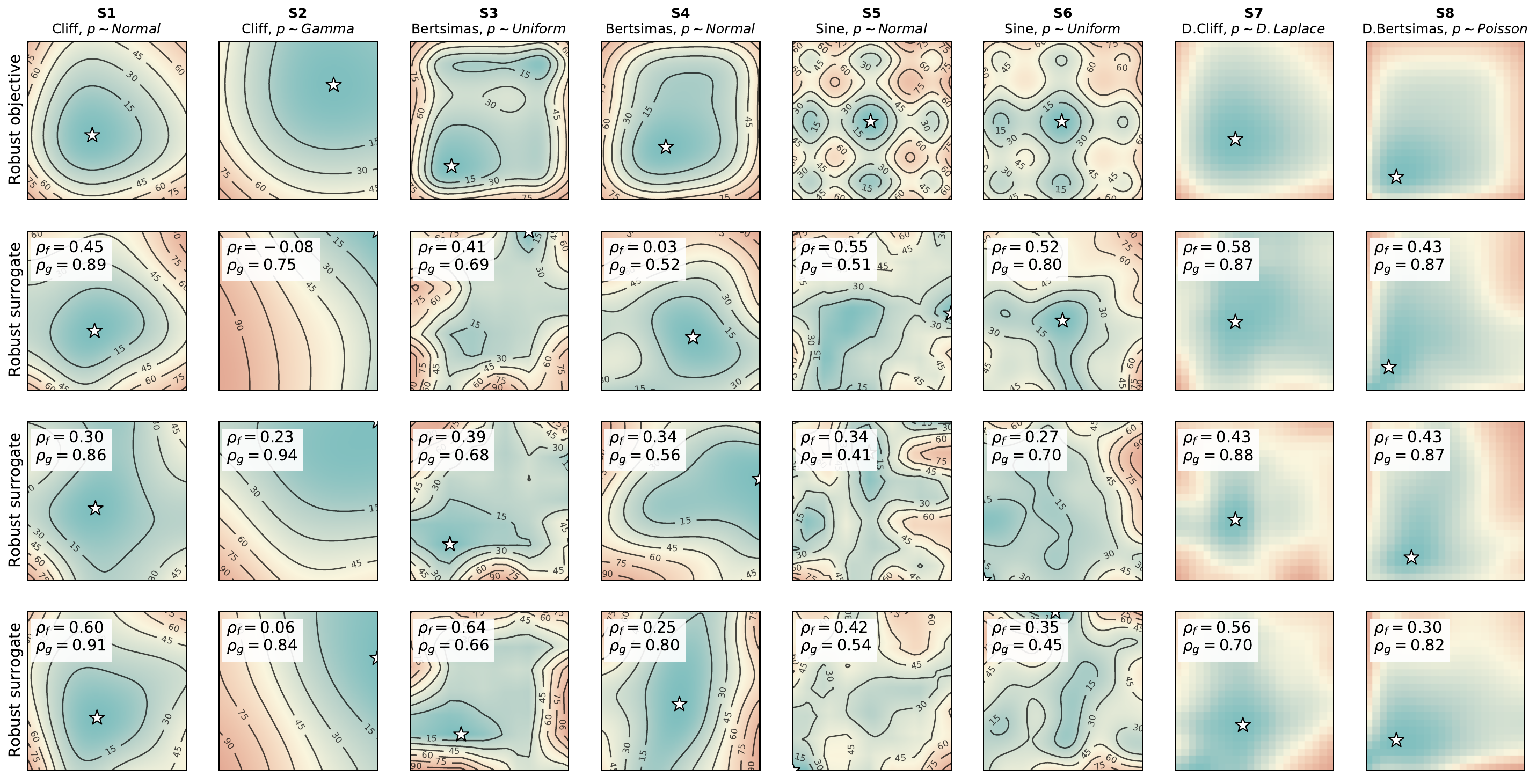}
    \caption{\golem's robustness estimates based on noisy data. The first row shows the true robust objective functions. The other rows show \golem estimates based on three different noisy datasets $\widetilde{\mathcal{D}}_K$. In all cases, $64$ datapoints were sampled under severe input noise, according to the uncertainties defined in Table~\ref{tab:surfaces}. On each plot, the Spearman's correlation between true and observed objective function values for these $64$ datapoints is reported as $\rho_f$. The correlation between \golem's robustness estimates (based on these noisy datapoints) and the true robust objective values is reported as $\rho_g$. In the vast majority of cases, $\rho_g$ was much larger than $\rho_f$, which means that \golem was able to recover correlations with the robust objective despite having to rely on poorly informative samples of the objective function. The most striking example was observed for surface S2, where in one case there was a negative correlation between the true and sampled objective function values ($\rho_f = -0.08$), while \golem's predictions showed a strong positive correlation ($\rho_g = 0.75$).}
    \label{si_fig:corr_fx_vs_gx}
\end{figure*}

As discussed in section~\ref{section:benchmarks_noisy}, when queries are noisy, building an accurate surrogate model is challenging because the objective function is not evaluated at the desired queried locations. As a consequence, the data $\widetilde{\mathcal{D}}_K = \{ \bm{x}_k, \widetilde{f}_k \} _{k=1}^K$ available to train the surrogate model is mismatched. However, while \golem does not take into account input noise at training time, it does so at the inference stage when estimating robustness (i.e., $g(\bm{x})$). \golem may be seen as trying to estimate where the value of $\widetilde{f}_k$ for the query $\bm{x}_k$ might have come from. Because of this, \golem is able to recover reasonably accurate estimates of $g(\bm{x})$ even when the correlation between true, $f_k$, and observed, $\widetilde{f}_k$, objective function values is lost. This effect is shown in Figure~\ref{si_fig:corr_fx_vs_gx}. The first row shows the true robust objective functions that we introduced in Figure~\ref{fig:grids}, while the other rows show three \golem estimates based on a dataset $\widetilde{\mathcal{D}}_K$ comprised of $64$ datapoint collected under severe input uncertainty (as defined in Table~\ref{tab:surfaces}). As a consequence, the correlation between true and observed objective function values, for the chosen queries locations is low, and in some cases effectively lost. However, \golem manages to recover moderate ($0.4-0.7$) to strong ($0.7-0.9$) correlations between its robustness estimates and the true robust objective values.

\subsection{Influence of the type and size of tree ensemble}
\label{section:si_forests}

\golem can be used with several tree ensemble algorithms. We tested the performance of \golem where the surrogate function is modeled with regression trees, random forest\cite{Breiman:2001}, and extremely randomized trees\cite{Geurts:2006}. The \textit{scikit-learn}\cite{scikit-learn} implementations of these algorithms were used (\texttt{DecisionTreeRegressor}, \texttt{RandomForestRegressor}, and \texttt{ExtraTreesRegressor}, respectively). In addition to testing the performance of a single regression tree, we also fitted ensembles of 10, 20, and 50 trees for all above-mentioned algorithms. Note that the regression tree algorithm used is not fully deterministic, such that different trees in the ensemble can correspond to different surrogate models. While the input dataset is not bootstrapped (like in random forest)\cite{Breiman:2001} and thresholds for splitting nodes are not chosen at random (like in extremely randomized trees)\cite{Geurts:2006}, multiple splits can provide the same mean-square-error improvement, and a specific split is then chosen at random among these.

\begin{figure*}[htb]
    \centering
    \includegraphics[width=1.0\textwidth]{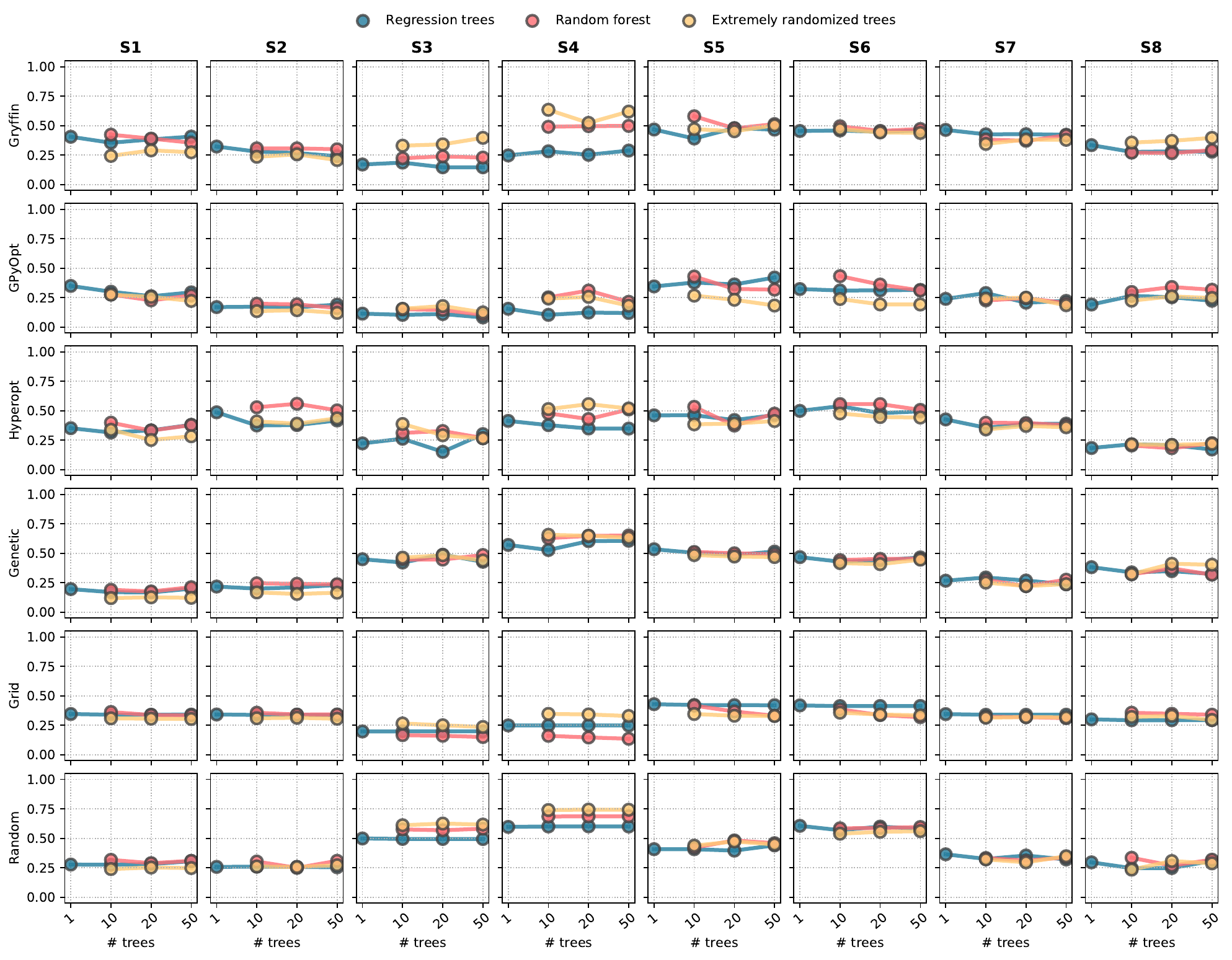}
    \caption{Influence of the type and size of tree ensemble on \golem's performance for optimizations in the noiseless query setting. Shown are cumulative regret values, normalized within each subplot, and averaged across 50 repeated optimization runs. Each subplot refers to optimizations performed on a different benchmark surface and with a different algorithm in conjunction with \golem.}
    \label{si_fig:tree_dependance1}
\end{figure*}

\begin{figure*}[htb]
    \centering
    \includegraphics[width=1.0\textwidth]{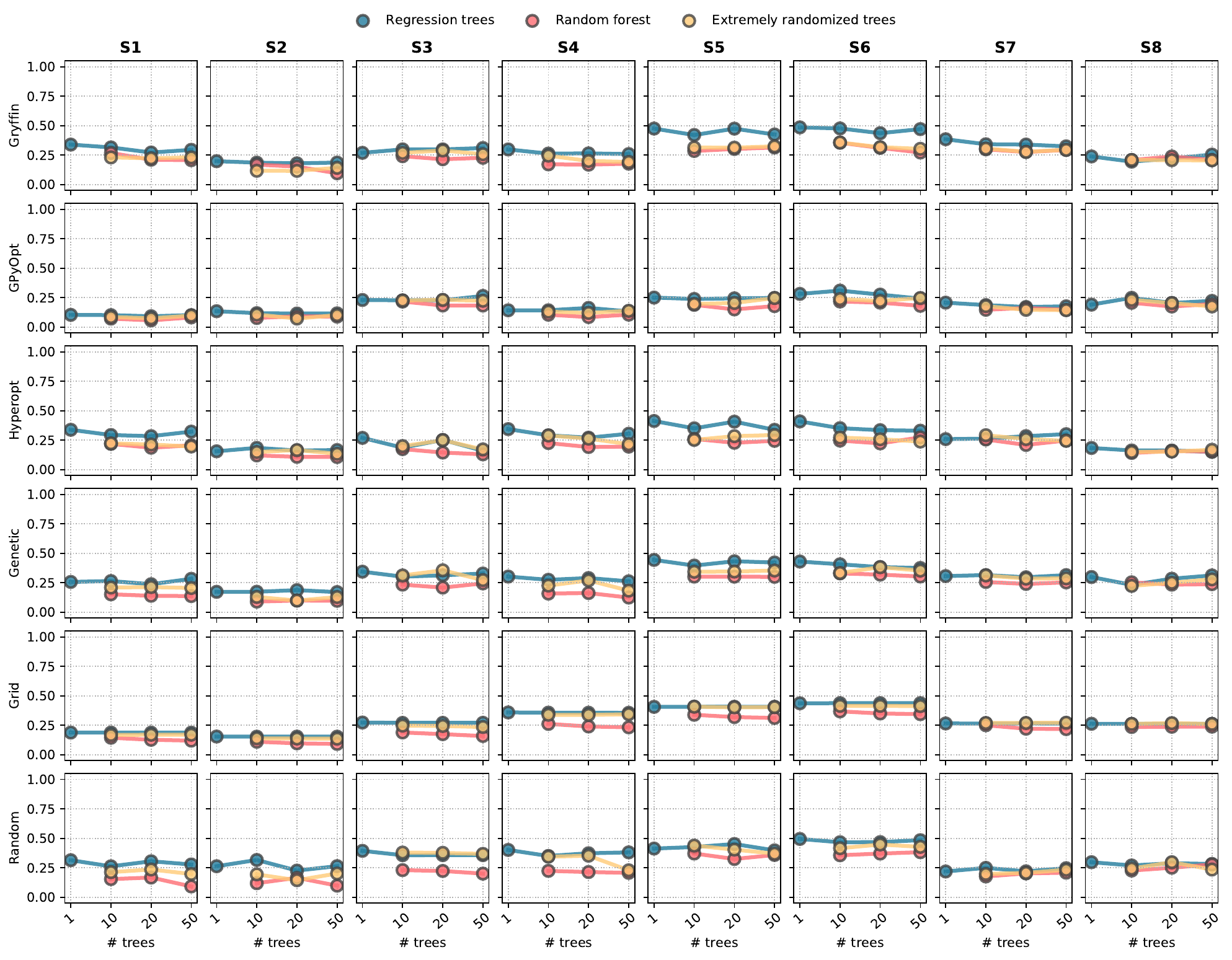}
    \caption{Influence of the type and size of tree ensemble on \golem's performance for optimizations in the noisy query setting. Shown are cumulative regret values, normalized within each subplot, and averaged across 50 repeated optimization runs. Each subplot refers to optimizations performed on a different benchmark surface and with a different algorithm in conjunction with \golem.}
    \label{si_fig:tree_dependance2}
\end{figure*}

Figures~\ref{si_fig:tree_dependance1} and \ref{si_fig:tree_dependance2} provide a summary of \golem's relative performance when using (i) different tree-based machine learning models (regression trees, random forest, and extremely randomized trees), and (ii) ensemble of trees of different sizes ($1$, $10$, $20$, $50$). Figures~\ref{si_fig:tree_dependance1} shows the normalized cumulative regret values (section \ref{section:si_regret}) for the results obtained with the six optimization algorithms tested, on the eight benchmark functions employed, and in the noiseless query setting. Figures~\ref{si_fig:tree_dependance2} shows the same results, but for optimizations in the noisy query setting. Note that, because of the normalization of the cumulative regrets, results can be compared only within, and not across, subplots.

\begin{figure*}[htb]
    \centering
    \includegraphics[width=1.0\textwidth]{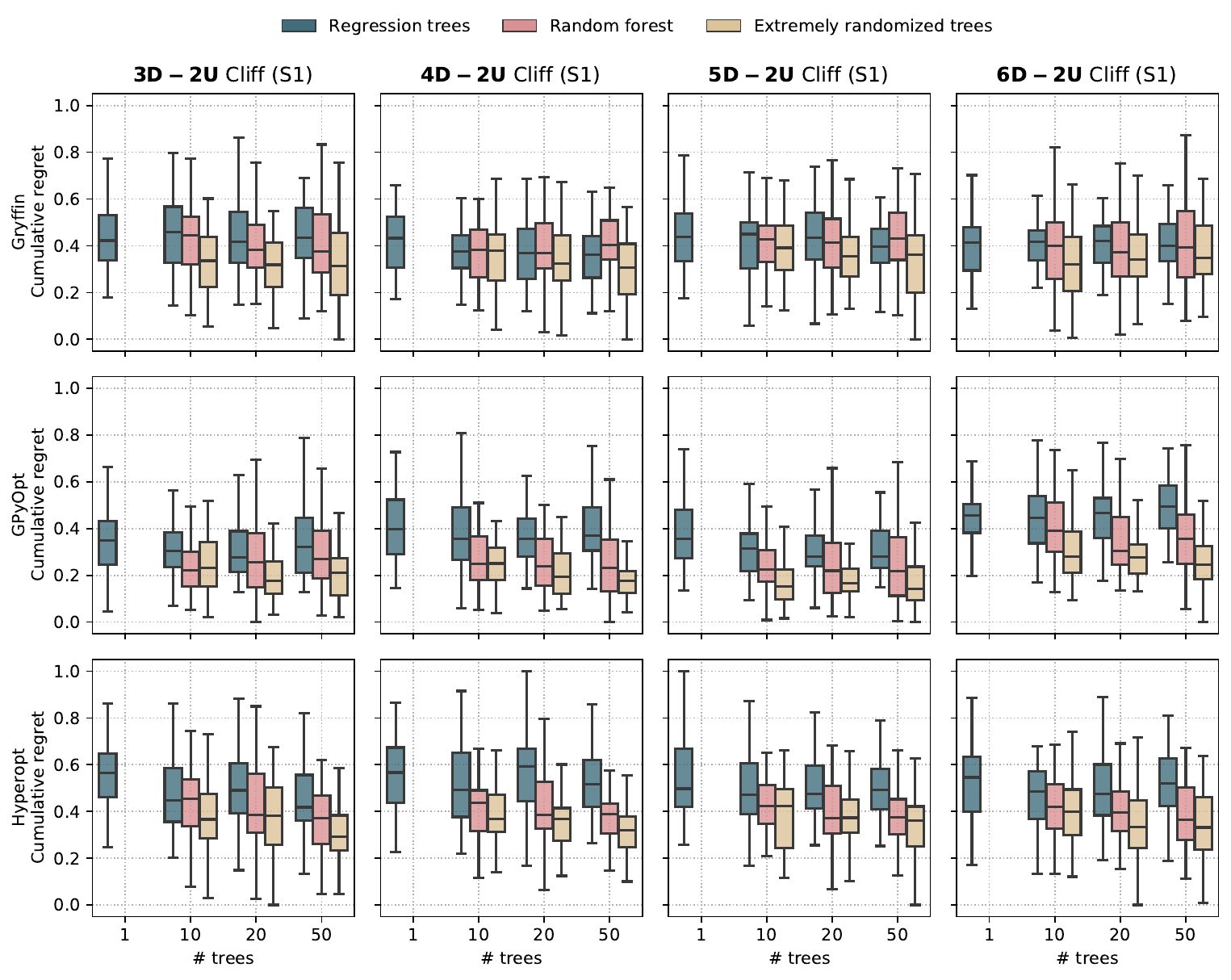}
    \caption{Influence of the type and size of tree ensemble on \golem's performance for high-dimensional optimizations in the noiseless query setting. The distributions of cumulative regret values, normalized within each subplot, across $50$ repeated optimization runs are shown. Each subplot refers to optimizations performed on surfaces of increasing dimensionality and with a different algorithm.}
   \label{si_fig:highdim_noiseless}
\end{figure*}

\begin{figure*}[htb]
    \centering
    \includegraphics[width=1.0\textwidth]{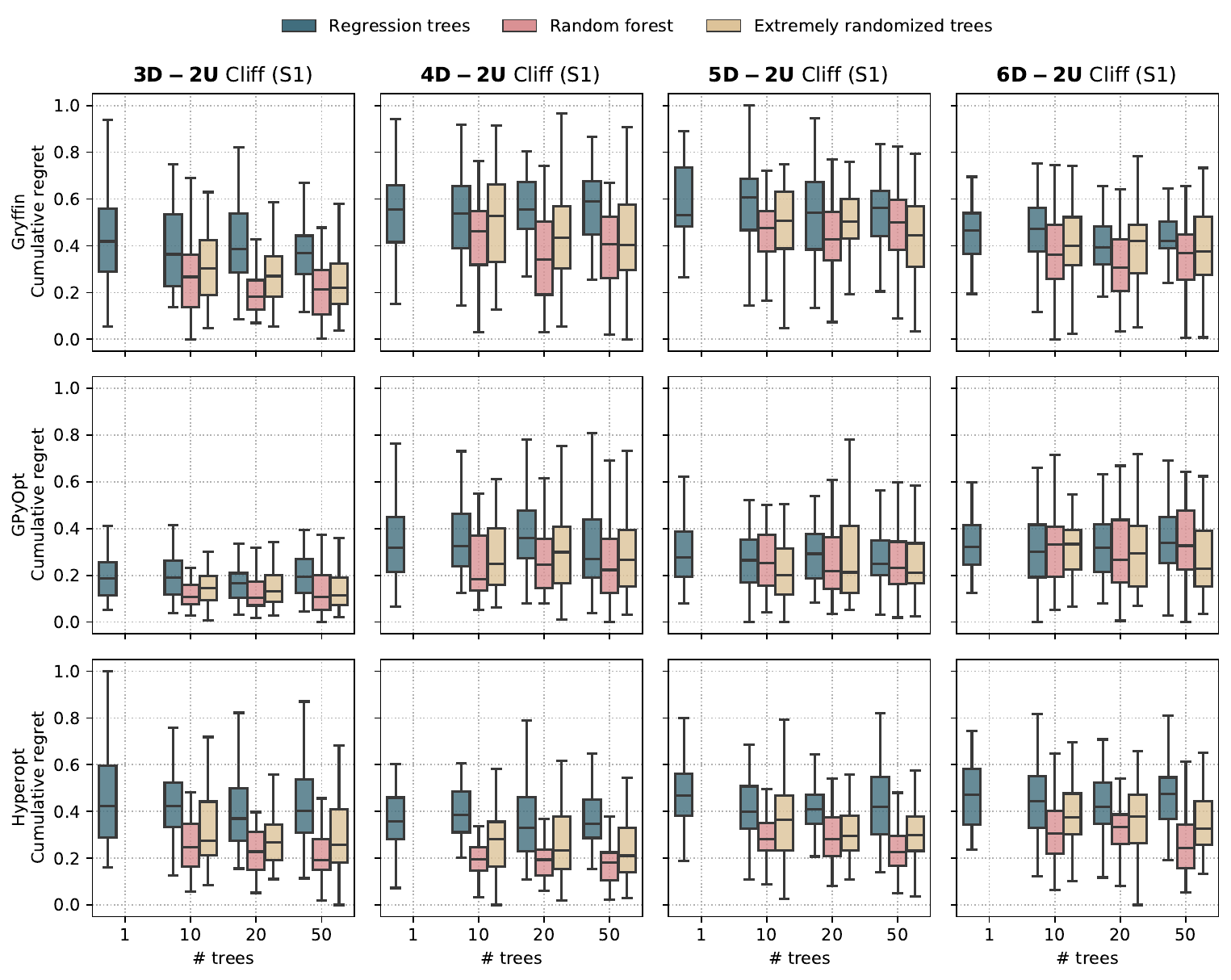}
    \caption{Influence of the type and size of tree ensemble on \golem's performance for high-dimensional optimizations in the noisy query setting. The distributions of cumulative regret values, normalized within each subplot, across $50$ repeated optimization runs are shown. Each subplot refers to optimizations performed on surfaces of increasing dimensionality and with a different algorithm.}
    \label{si_fig:highdim_noisy}
\end{figure*}

Figures \ref{si_fig:highdim_noiseless} and \ref{si_fig:highdim_noisy} provide a summary of \golem's relative performance when using different tree-based models of different size, on high-dimensional benchmark surfaces. All these surfaces are higher-dimensional versions of the surface S1 (from three to six dimensions), where the first two dimensions are uncertain, while the additional ones are always considered to be noiseless. In these high-dimensional tests, it is possible to notice how surrogate models based on random forest and extremely randomized trees provided slightly better performance than regression trees.

\clearpage

\clearpage

\clearpage

\subsection{Influence of the number of uncertain variables}
\label{section:si_num_dims}
\begin{figure*}[htb]
    \centering
    \includegraphics[width=0.95\textwidth]{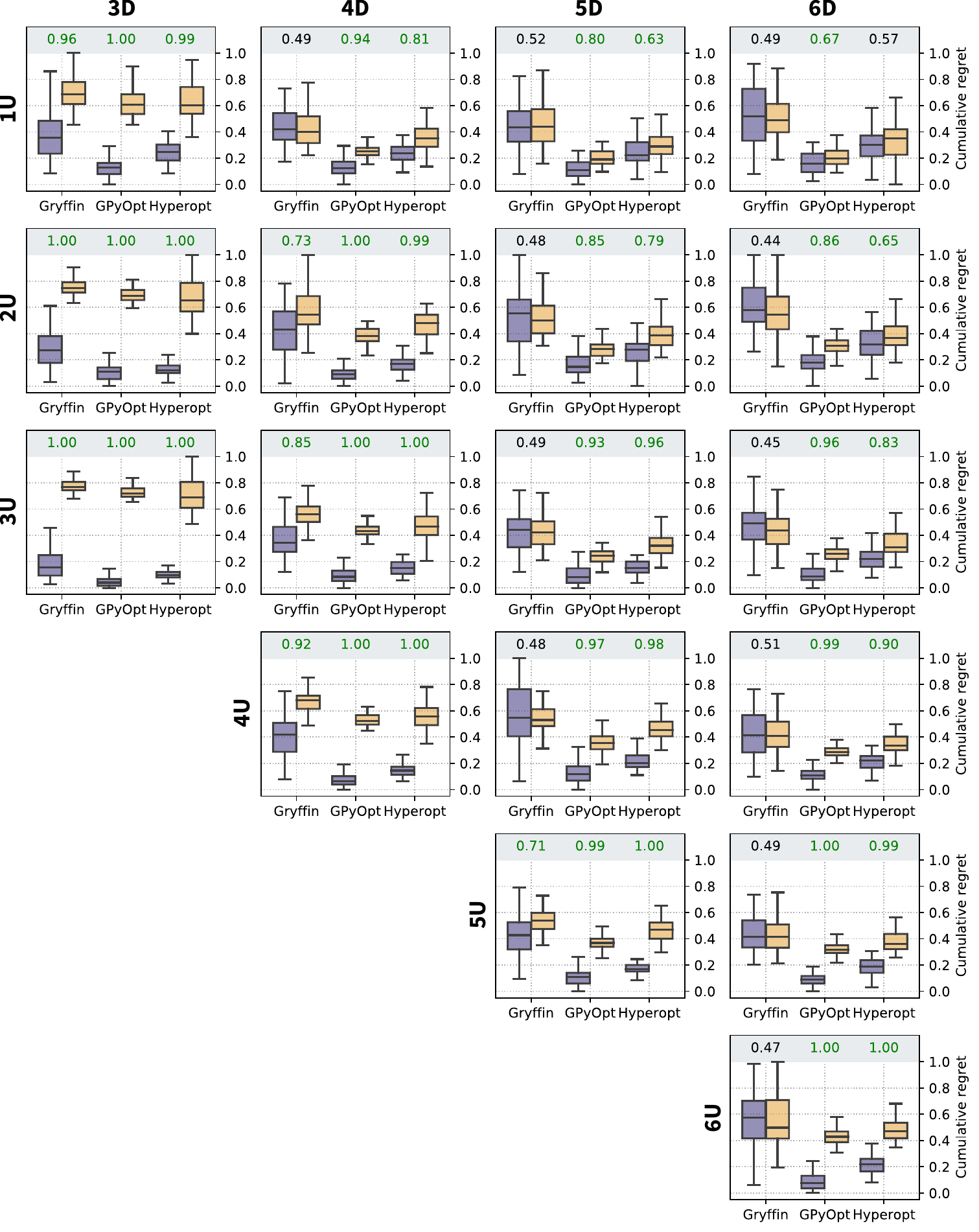}
    \caption{Relative comparison of optimization performance obtained with and without \golem on the surface S1 with varying dimensions (3D$-$6D) and number of uncertain inputs (1U$-$6U) in the noiseless query setting. The regret distributions shown were obtained from optimizations that used \golem with an ensemble of $50$ extremely randomized trees as the surrogate model. The boxes show the first, second, and third quartiles of the data, with whiskers extending up to 1.5 times the interquartile range. Results obtained with \golem are shown in purple, and those obtained without \golem in yellow. The probability of obtaining better performance with \golem, with the algorithms tested, is reported above each box. Statistically significant results ($\alpha=0.05$) are highlighted in green (significant improvement when using \golem) and red (significant deterioration when using \golem).}
    \label{si_fig:highdim_highunc_noiseless}
\end{figure*}

\begin{figure*}[htb]
    \centering
    \includegraphics[width=0.95\textwidth]{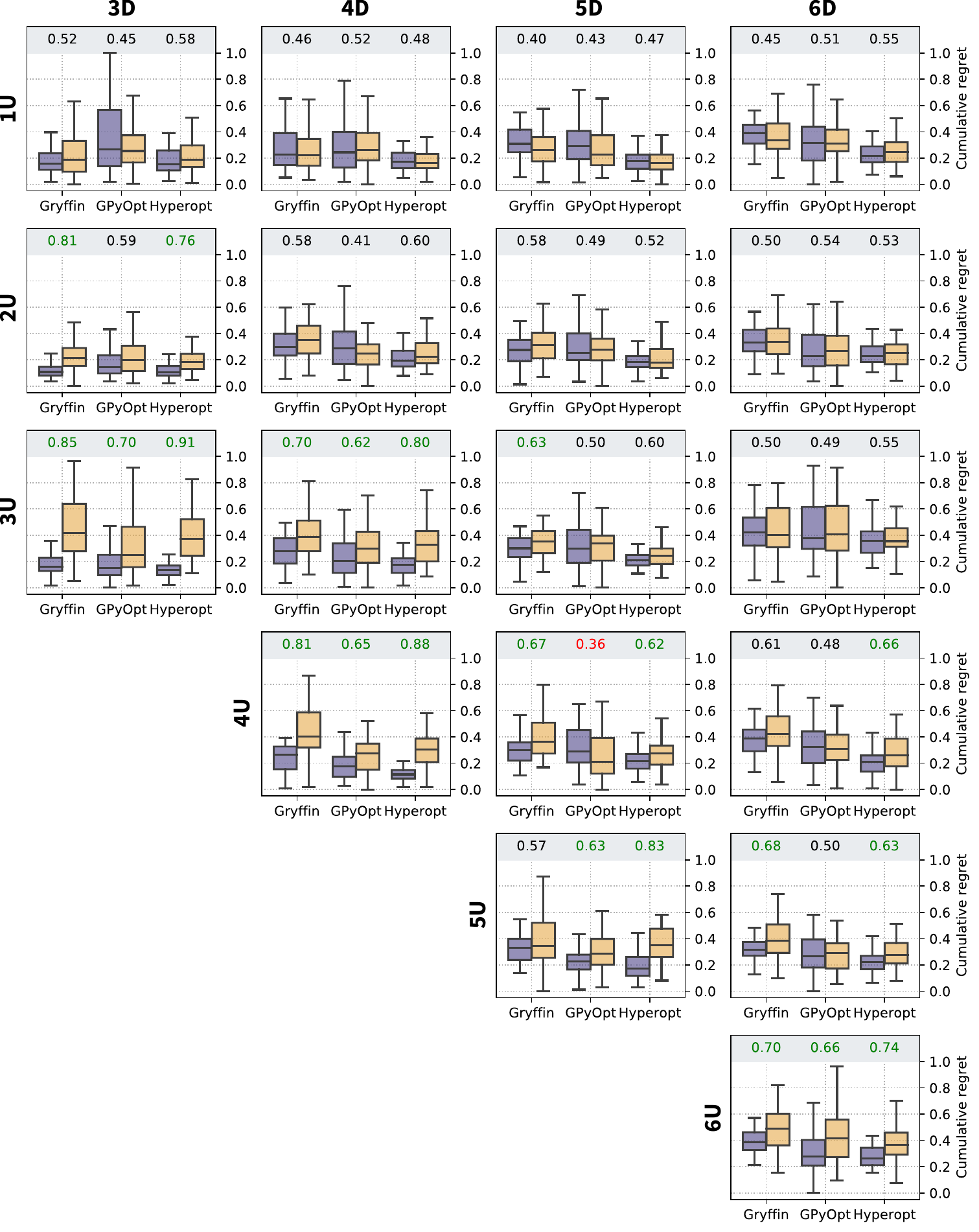}
    \caption{Relative comparison of optimization performance obtained with and without \golem on the surface S1 with varying dimensions (3D$-$6D) and number of uncertain inputs (1U$-$6U) in the noisy query setting. The regret distributions shown were obtained from optimizations that used \golem with an ensemble of $50$ extremely randomized trees as the surrogate model. The boxes show the first, second, and third quartiles of the data, with whiskers extending up to 1.5 times the interquartile range. Results obtained with \golem are shown in purple, and those obtained without \golem in yellow. The probability of obtaining better performance with \golem, with the algorithms tested, is reported above each box. Statistically significant results ($\alpha=0.05$) are highlighted in green (significant improvement when using \golem) and red (significant deterioration when using \golem).}
    \label{si_fig:highdim_highunc_noisy}
\end{figure*}

The larger the dimensionality of the problem, and the larger the number of uncertain input variables, the more challenging the robust optimization task is. We tested how \golem's performance is affected by the presence of additional noise-free and noisy input variables. All surfaces used in these tests are higher-dimensional versions of the surface S1 (from three to six dimensions), where between one and all of the available input variables are noisy. Figures~\ref{si_fig:highdim_highunc_noiseless} and \ref{si_fig:highdim_highunc_noisy} show the normalized cumulative regrets obtained for optimizations in the noiseless and noisy query setting, respectively. These results are discussed in \ref{section:benchmarks_highdim}. In summary, we find that \golem is effective also on higher-dimensional surfaces. In fact, the benefits of using \golem become more marked the higher the number of uncertain inputs present in the optimization domain (i.e., the more uncertainty being present overall). On the other hand, given a fixed number of uncertain input variables, additional noiseless variables make it harder for \golem to enhance the performance of the optimization algorithm used. Importantly, \golem was almost never (one out of $108$ tests) found to be detrimental to optimization performance.

\clearpage

\clearpage

\section{Analysis and optimization of an HPLC protocol}
\label{section:si_hplc}

In this section we provide more details on the setup and results concerning the example application on the calibration of an HPLC protocol. In this example, the experimental HPLC response depends on six tunable parameters. These controllable input parameters (shown in Figure \ref{fig:hplc1}a) are the following: (P1) volume of the sample loop and internal volume of the 2-way 6-port valve; (P2) volume required to draw the sample to the 2-way 6-port valve; (P3) volume required to drive the sample plug from the sample loop, through the in-line mixer, and to the second valve; (P4) draw rate of the sample pump; (P5) push rate of the push pump; (P6) time waited after drawing sample and before switching the first selection valve (to allow for equilibration of cavitation bubbles in the sample line and syringe). As discussed in the main text, \golem may be used to retrospectively analyze the experimental results, or to optimize the protocol assuming no prior knowledge.

\subsection{Interaction between input uncertainties and optimum location}
\label{section:si_hplc_interac}

\begin{figure*}[htb]
    \centering
    \includegraphics[width=1.0\textwidth]{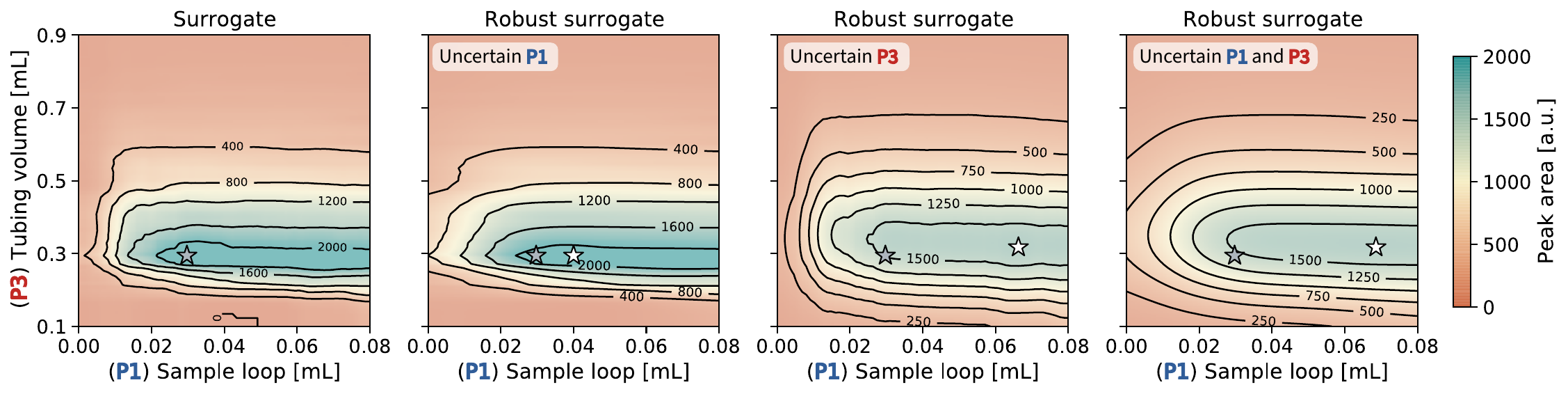}
    \caption{Effect of uncertainty in P1 and P3 on the optimum location of the HPLC protocol. \golem's surrogate models are shown against the input parameters P1 and P3, while the other parameters are fixed. The plot on the left-hand side shows \golem's surrogate model, while the other plots show the robust counterpart when assuming uncertainty in P1, P3, and P1 and P3. We assume a normally distributed uncertainty with standard deviation corresponding to $10$\% of the parameter range ($0.008$ for P1 and $0.08$ for P3). The distribution is furthermore truncated at zero to avoid non-physical values of P1 or P3. The location of the non-robust optimum is indicated by a gray star (as found by the experimental sample collected with highest peak area), while the location of the robust optima is indicated by a white star (as computed with \golem). These results show how uncertainty in P3 results in a large shift in optimum location along P1.}
    \label{si_fig:hplc_interac}
\end{figure*}

The uncertainty present in one input parameter affects the robustness merit of the solutions across the whole search space. As such, uncertainty in one parameter might affect the optimal setting for other parameters too. Figure \ref{si_fig:hplc_interac} shows such an example to visually clarify this statement. Based on a surrogate model built on $1386$ HPLC experiments, we can use \golem to investigate how the response surface is affected by uncertainty in the parameters P1 and P3. For ease of visualization, Figure \ref{si_fig:hplc_interac} shows surfaces only with respect to P1 and P3, with the other parameters fixed according to the best performing sample collected (P2 $\approx 0.03$ mL, P4 $\approx 2.4$ mL/min, P5 $\approx 107$ Hz, P6 $\approx 6.2$ s). We consider the presence of input uncertainty in P1 and P3 individually, and then in P1 and P3 together. In all cases, we assume a normally distributed uncertainty with standard deviation corresponding to $10$\% of the parameter range (i.e., $0.008$ for P1 and $0.08$ for P3). The distribution is furthermore truncated at zero to avoid non-physical values of P1 or P3. Assuming uncertainty only in P1, the location of the optimum is shifted to slightly higher values of P1 (Figure \ref{si_fig:hplc_interac}). However, when assuming uncertainty in P3, the location of the optimum is shifted considerably towards higher values of P1, while leaving almost unaffected the location with respect to P3, the uncertain parameter. This effect is due to objective function dropping slightly more steeply towards zero at low P3 values also when P1 is low, while having a slightly broader maxima in the P3 dimension for higher values of P1. The net result of this effect is that, when considering uncertainty in both P1 and P3, the location of the optimum is shifted primarily in P1, yet it is determined mainly by the uncertainty in P3. These types of interactions between variables are difficult to discover by simple visual inspection of the surrogate model, and is one of the tasks in which the use of \golem proves useful.

\subsection{Optimization of a noisy HPLC protocol}
\label{section:si_hplc_interac}

In this example application, we assumed the presence of noise in parameters P1 and P3 while attempting to optimize the HPLC sampling protocol. We assumed this noise to be normally distributed and truncated at zero, with standard deviation of $0.008$ mL for P1, and $0.08$ mL for P3. The HPLC experiments were simulated with \textsc{Olympus}\cite{Hase:2020_olympus}, which emulates the experimental HPLC response based on its six tunable parameters via a Bayesian Neural Network. The goal of the optimization was to achieve a protocol returning an expected peak area, $\mathbb{E}[Area]$, of at least $1000$ a.u. As a secondary objective, we wanted to minimize the output variability, $\sigma[Area]$, as much as possible, as long as $\mathbb{E}[Area] > 1000$ a.u. \golem was used to estimate both the $\mathbb{E}[Area]$ and $\sigma[Area]$ during the optimization (Figure \ref{fig:hplc-opt}a), using $200$ extremely randomized trees\cite{Geurts:2006} as the surrogate model. The \textit{Chimera}\cite{Hase:2018_chimera} scalarizing function was used to create a robust, multi-objective function to be optimized.

\begin{figure*}[htb]
    \centering
    \includegraphics[width=0.8\textwidth]{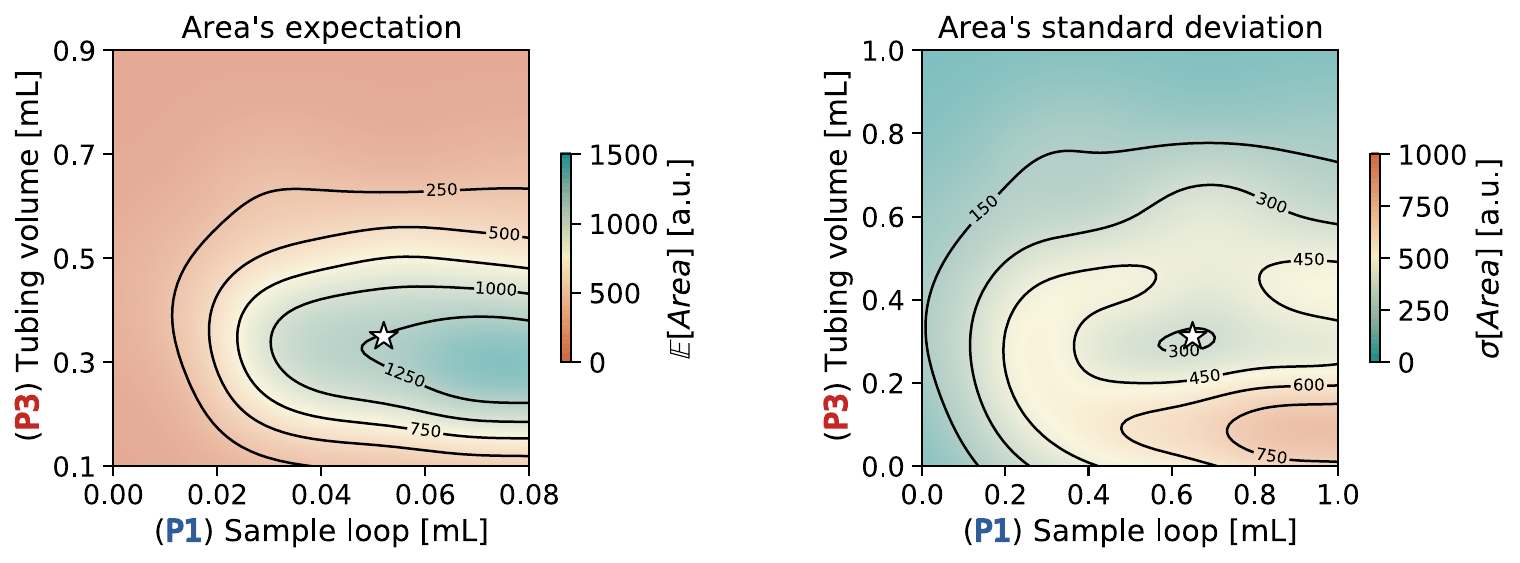}
    \caption{Location of the true robust optimum identified and behavior of $\mathbb{E}[Area]$ and $\sigma[Area]$ around this optimum. The location of the global optimum is marked by a white star. It is located at P1 $\approx 0.052$ mL, P2 $\approx 0.012$ mL, P3 $\approx 0.35$ mL, P4 $\approx 2.20$ mL/min, P5 $\approx 84$ Hz, P6 $\approx 5.9$ s, where $\mathbb{E}[Area] = 1256$ and $\sigma[Area] = 288$.}
    \label{si_fig:hplc_true_optimum}
\end{figure*}

Similar to what we did to obtain a ground truth for the robust objectives for the analytical surfaces (Section \ref{section:si_surfaces}), a close numerical approximation of $\mathbb{E}[Area]$ and $\sigma[Area]$ was obtained by using a dense grid of uniformly distributed samples across the optimization domain. In this case, we sampled $8^6 = 262,144$ points from the \textsc{Olympus} experiment emulator to build a reference \golem model. These samples were extended beyond the optimization domain in P1 and P3 by two standard deviations. The approximate location of true robust optimum (Figure \ref{si_fig:hplc_true_optimum}) was found with \textit{Hyperopt} by optimizing the true robust, multi-objective function directly over $1000$ iterations.

\begin{figure*}[htb]
    \centering
    \includegraphics[width=1\textwidth]{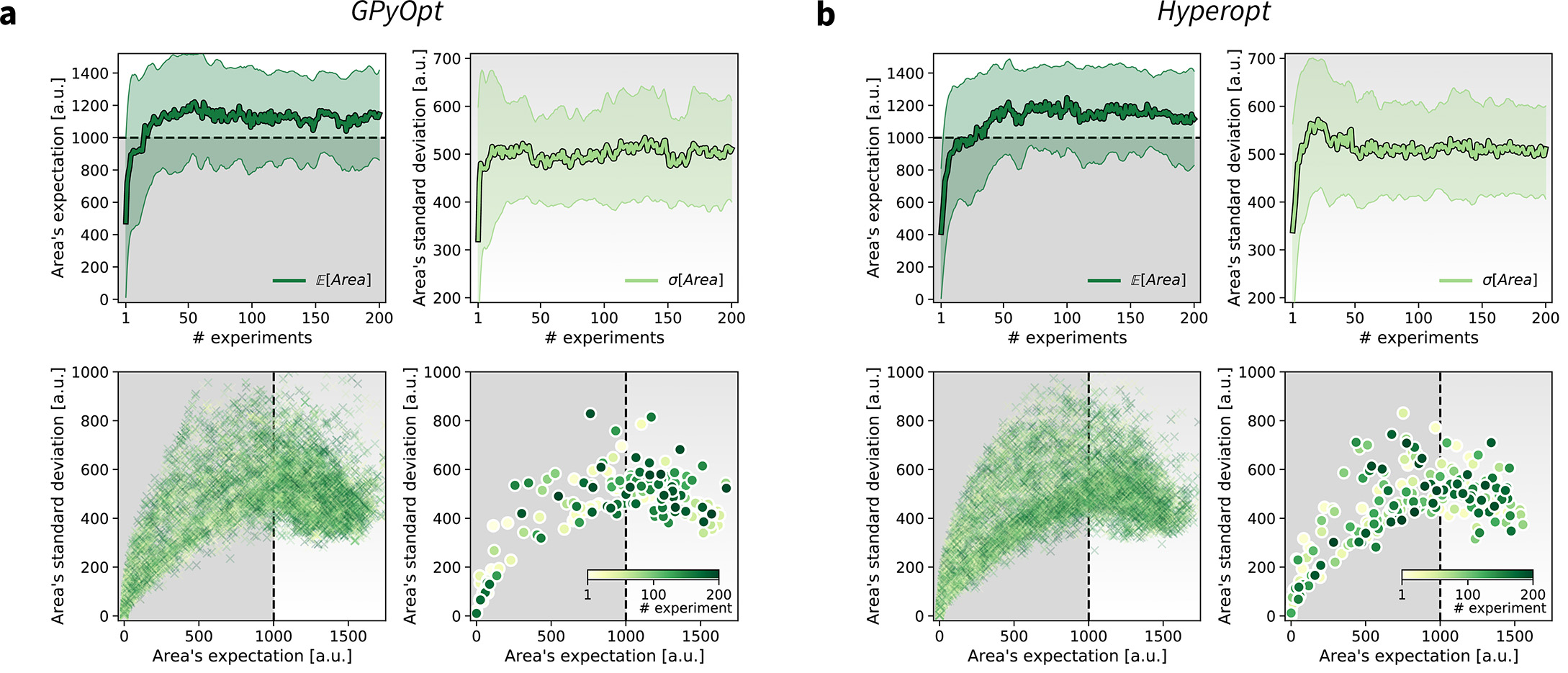}
    \caption{Results of 50 optimization repeats performed with (a) \textit{GPyOpt} and (b) \textit{Hyperopt}. In both cases, optimization traces for the primary and secondary objectives are shown (average and standard deviation). All objective function values sampled during the optimization runs are shown in the bottom-left panels. Objective function values sampled during an example optimization run are shown in the bottom-right panels, with each experiment color-coded (yellow to dark green) to indicate at which stage of the optimization it was performed.}
    \label{si_fig:hplc_opt_gpyopt_hopt}
\end{figure*}

\begin{figure*}[htb]
    \centering
    \includegraphics[width=0.9\textwidth]{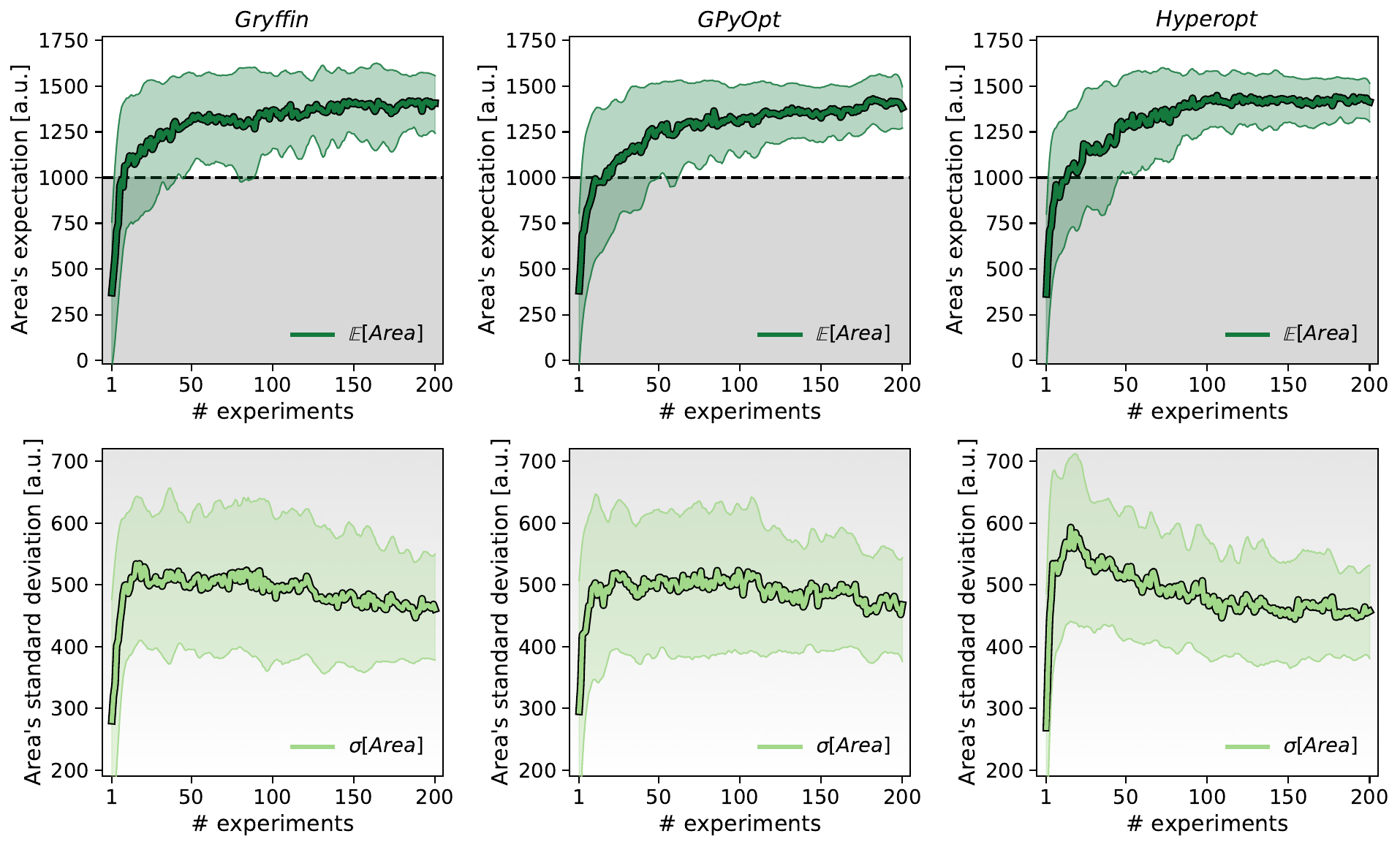}
    \caption{Traces of 50 optimization repeats in which the primary objective was constrained to a lower-bound estimate of the peak area's expectation, $\mathbb{E}[Area] - 1.96 \times \sigma(\mathbb{E}[Area])$, corresponding to the lower bound of the 95\% confidence interval of \golem's estimates. After $200$ experiments, \textit{Gryffin}, \textit{GPyOpt}, and \textit{Hyperopt} all correctly identified solutions with $\mathbb{E}[Area] > 1000$ a.u. in all $50$ repeated optimization runs.}
    \label{si_fig:hplc_opt_lcb}
\end{figure*}

\begin{figure*}[htb]
    \centering
    \includegraphics[width=0.9\textwidth]{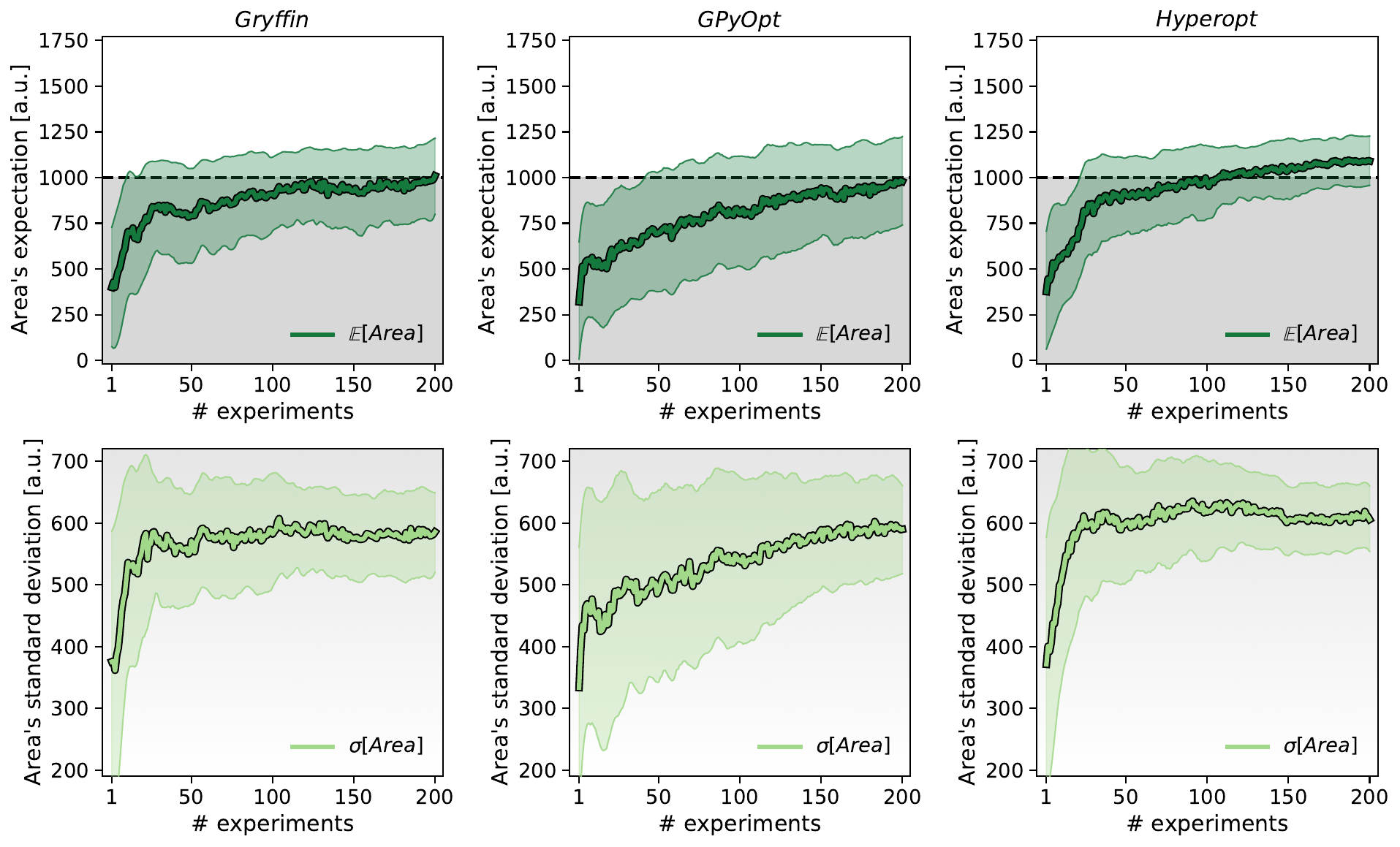}
    \caption{Traces of 50 optimization repeats in which all six input parameters were noisy. The primary objective was constrained to a lower-bound estimate of the peak area's expectation, $\mathbb{E}[Area] - 1.96 \times \sigma(\mathbb{E}[Area])$, corresponding to the lower bound of the 95\% confidence interval of \golem's estimates. With more overall noise in the input experimental conditions, the optimization takes longer than in the previous example with only two noisy variables. However, with \golem, all approaches still managed to optimize the peak area's expectation (i.e., $\mathbb{E}[Area]$ increases with more experiments performed and, on average, reaches the targeted value of $1000$ a.u.). After $200$ experiments, \textit{Gryffin} correctly identified solutions with $\mathbb{E}[Area] > 1000$ a.u. in $42$\% of the optimization runs, \textit{GPyOpt} in $70$\% of the optimization runs, and  \textit{Hyperopt} in $78$\%.}
    \label{si_fig:hplc_opt_lcb_all_noisy}
\end{figure*}

	\putbib[main]
\end{bibunit}

\end{document}